\documentclass[11pt,twoside]{article}%
\usepackage{enumerate}
\usepackage{enumitem}
\usepackage{amsmath,amssymb,amsfonts}
\usepackage{amsfonts}%
\usepackage{graphicx}
\usepackage{color}
\usepackage{pdfsync}
\providecommand{\U}[1]{\protect\rule{.1in}{.1in}}
\newtheorem{thm}{Theorem}[section]
\newtheorem{lem}[thm]{Lemma}

\newtheorem{prop}[thm]{Proposition}
\newtheorem{corol}[thm]{Corollary}

\newenvironment{pf}[1][\bfseries Proof]{\noindent{#1.} }{\hfill \rule{0.5em}{0.5em}\\}
\newcommand{\fim}{\hfill\rule{2mm}{2mm}}
\numberwithin{equation}{section}
\begin{document}
\title{Multiplicity and concentration behavior of solutions for a quasilinear problem involving $N$-functions via penalization method}

\author{Claudianor O. Alves\thanks{C.O. Alves was partially supported by CNPq/Brazil  304036/2013-7  and INCT-MAT, \hspace*{.7cm} e-mail:  coalves@dme.ufcg.edu.br},~~  Ailton R. Silva\thanks{A.R. Silva,~ e-mail: ailton@dme.ufcg.edu.br}\\
Universidade Federal de Campina Grande\\
 Unidade Acad\^emica de Matem\'atica - UAMat\\
 CEP: 58.429-900 - Campina Grande - PB - Brazil}
\date{}
\maketitle

\begin{abstract}
In this work, we study the existence, multiplicity and concentration of positive solutions for the following class of quasilinear problem:
\[
- \Delta_{\Phi}u + V(\epsilon x)\phi(\vert u\vert)u  =  f(u)\quad \mbox{in} \quad \mathbb{R}^{N},
\]
where $\Phi(t) = \int_{0}^{\vert t\vert}\phi(s)sds$ is a N-function, $ \Delta_{\Phi}$ is the $\Phi$-Laplacian operator, $\epsilon$ is a positive parameter, $ N\geq 2$, $V : \mathbb{R}^{N} \rightarrow \mathbb{R} $ is a continuous function and $f : \mathbb{R} \rightarrow \mathbb{R} $ is a $C^{1}$-function.

\end{abstract}

{\scriptsize \textbf{2000 Mathematics Subject Classification:} 35A15, 35B40, 35J62, 46E30 }

{\scriptsize \textbf{Keywords:} Variational Methods, Behavior of Solutions, Quasilinear Problems, Orlicz-Sobolev Space}

\section{Introduction}

Many recent studies have focused on the nonlinear Schr\"{o}dinger
equation
$$
i\epsilon \displaystyle \frac{\partial \Psi}{\partial t}=-\epsilon^{2}\Delta
\Psi+(V(z)+E)\Psi-f(\Psi)\,\,\, \mbox{for all}\,\,\, z \in
\mathbb{R}^{N},\eqno{(NLS)}
$$
where $N \geq 1$, $\epsilon > 0$ is a parameter and $V,f$ are continuous function verifying some conditions. This class of equation is one of the main objects of the quantum physics, because it appears in problems involving nonlinear optics, plasma physics and condensed matter physics.

Knowledge of the solutions for the elliptic equation
$$
\ \  \left\{
\begin{array}{l}
- \epsilon^{2} \Delta{u} + V(z)u=f(u)
\ \ \mbox{in} \ \ \mathbb{R}^{N},
\\
u \in H^{1}(\mathbb{R}^{N}),
\end{array}
\right.
\eqno{(S)_{\epsilon}}
$$
or equivalently
$$
\ \  \left\{
\begin{array}{l}
-\Delta{u} + V(\epsilon z)u=f(u)
\ \ \mbox{in} \ \ \mathbb{R}^{N},
\\
u \in H^{1}(\mathbb{R}^{N}),
\end{array}
\right.
\eqno{(S')_{\epsilon}},
$$
has a great importance in the study of standing-wave solutions
of $(NLS)$. In recent years, the existence and concentration of
positive solutions for general semilinear elliptic equations
$(S)_\epsilon$ have been extensively studied, see for example, Floer and Weinstein \cite{FW}, Oh
\cite{O1,O2}, Rabinowitz \cite{rabinowitz}, Wang \cite{WX}, Ambrosetti and Malchiodi \cite{AM}, Ambrosetti, Badiale and Cingolani  \cite{ABC}, Floer and Weinstein \cite{FN}, del Pino and Felmer \cite{DF1}  and their references.

In the above mentioned papers, the existence,  multiplicity and concentration of positive solutions have been obtained in connection with the geometry of the function $V$. In some of them, it was proved that the maximum points of the solutions are close to the set
$$
\mathcal{V} = \left\{x \in \mathbb{R}^{N} \ : \ V(x) = \min_{z \in \mathbb{R}^{N}}V(z) \right\},
$$
when $\epsilon$ is small enough. Moreover, in a lot of problems, the multiplicity of solution is related to topology richness of $\mathcal{V}$.

In \cite{rabinowitz}, by a mountain pass argument,  Rabinowitz proves the
existence of positive solutions of $(S)_{\epsilon}$, for $\epsilon > 0$
small, whenever
$$
\liminf_{|z| \rightarrow \infty} V(z) > \inf_{z \in
	\mathbb{R}^N}V(z)=V_{0} >0. \eqno{(R)}
$$
Later Wang \cite{WX} showed that these solutions concentrate at global minimum points of $V$ as  $\epsilon$ tends to 0.

In \cite{DF1}, del Pino and Felmer have found solutions which concentrate around local minimum of $V$ by introducing a penalization method.
More precisely, they assume that
$$
V(x)\geq \inf_{z \in \mathbb{R}^N}V(z)=V_{0}>0 \,\,\, \mbox{for all}\,\, \ x \in \mathbb{R}^{N}\eqno{(V_{0})}
$$
and there is an open and bounded set $\Omega \subset \mathbb{R}^{^N}$ satisfying
$$
\inf_{z \in \Omega}V(z)< \min_{z \in
	\partial \Omega}V(z). \eqno{(V_{1})}
$$

The existence, multiplicity  and concentration of positive solution have been also considered for quasilinear problems of the type
$$
-\Delta_p{u}+V(\epsilon x)|u|^{p-2}u=f(u), \,\,\, \mbox{in} \,\,\, \mathbb{R}^{N}
$$
and
$$
-\Delta_p{u}-\Delta_q{u}+V(\epsilon x)(|u|^{p-2}u + |u|^{q-2}u)=f(u), \,\,\, \mbox{in} \,\,\, \mathbb{R}^{N}.
$$
Related to this subject, we cite the papers Alves and Figueiredo \cite{AF2005,AG2006}, Benouhiba and Belyacine \cite{BB}, Cammaroto and Vilasi \cite{CV}, Chaves, Ercole and Miyagaki \cite{CEM}, Figueiredo \cite{F1}, Li and Liang \cite{LL} and their references.

Recently, in \cite{AS}, Alves and Silva have showed the existence,  multiplicity and concentration of positive solutions for the following class of quasilinear problems
\begin{align}
\left\{
\begin{array}
[c]{rcl}%
- \Delta_{\Phi}u + V(\epsilon x)\phi(\vert u\vert)u & = & f(u)~  \mbox{in}~ \mathbb{R}^{N},\\
u \in W^{1, \Phi}(\mathbb{R}^{N}), &  &
\end{array}
\right. \tag{$ P_{\epsilon} $}\label{P2}%
\end{align}
where $ N\geq 2$, $\epsilon$  is positive parameter, the operator $ \Delta_{\Phi} u  = \mathrm{div}(\phi(\vert\nabla u\vert)\nabla u) $,  where $\Phi(t) = \int_{0}^{\vert t\vert}\phi(s)sds$, named $\Phi$-Laplacian, is a natural extension of the $p$-Laplace operator, with $ p $ being a positive constant and $V:\mathbb{R}^{N} \to \mathbb{R}$ is a continuous function verifying $(R)$.

This type of operator arises in a lot of applications, such as \\
\noindent {\it Nonlinear Elasticity:} \, $ \Phi(t) = (1+t^{2})^{\alpha}-1, \alpha \in (1, \frac{N}{N-2})$, \\
\noindent { \it Plasticity:} \, $ \Phi(t) = t^{p}\ln(1+t), 1< \frac{-1+\sqrt{1+4N}}{2}<p<N-1, N\geq 3$, \\
\noindent {\it Non-Newtonian Fluid:  } \, $\Phi(t) = \frac{1}{p}|t|^{p}$ for $p>1$,\\
\noindent { \it Plasma Physics: } \, $\Phi(t) = \frac{1}{p}|t|^{p} + \frac{1}{q}|t|^{q}$ where $1<p<q<N$ with $q \in (p, p^{*}).$ \\

The reader can find more details involving this subject in \cite{Db}, \cite{F1}, \cite{FN2} and their references.

Actually, we have observed that there are interesting papers studying the existence of solution for (\ref{P2}) when $\epsilon=1$, we would like to cite the papers \cite{AGJ}, \cite{BBR2}, \cite{FN}, \cite{FN1}, \cite{FN2}, \cite{MR1}, \cite{MR2}, \cite{J} and references therein. However, the authors know only the paper \cite{AS}, where the existence, multiplicity and concentration of solution has been considered for a $\Phi$-Laplacian equation.

Motivated by \cite{AF2005}, \cite{AS} and \cite{DF1}, in the present paper we study the existence,  multiplicity and concentration of solution for $(P_\epsilon)$, by supposing that $V$ verifies the conditions $(V_{0})$-$(V_{1})$.

In this work, we complete the study made \cite{AS}, because we are considering another geometry on $V$. Related to \cite{AF2005, DF1},  we enlarge their study, in the sense that, we obtain the same type of results for a large class of operators. In the proof of our results, we will work with N-function theory and Orlicz-Sobolev spaces. Since we are working with a general class of function $\Phi$, some estimates explored in \cite{AF2005} and \cite{DF1} cannot be repeated, then it was necessary to developed new estimates. For example, in \cite{AF2005}, it was used interaction Moser techniques, which does not work well in our case. To we overcome this difficulty, we adapt some arguments found in \cite{AF}, \cite{FUCHS}, \cite{Fusco}, \cite{LU} and \cite{ZF}. Here, we also modify the nonlinearity like \cite{DF1}, however the deformation is more technical, see Section 2 for details.

\vspace{0.5 cm}

Next, we will write our assumptions on $\phi$ and $f$.

\begin{flushleft}
\textbf{Conditions on $\phi$:}
\end{flushleft}
The function $\phi: [0, +\infty)\rightarrow [0, +\infty)$ is a $C^{1}$- function satisfying
\begin{enumerate}[label={($\phi_\arabic{*}$})]
\setcounter{enumi}{0}
\item\label{H1} $\phi(t)$, $(\phi(t)t)^{'}>0$, $t>0$.

\item\label{H2} There exist $l,m \in (1, N)$ such that
  \[
  l\leq m <l^{*}=\frac{Nl}{N-l}
  \]
  and
  \[
  l\leq \frac{\phi(t)t^{2}}{\Phi(t)}\leq m,  \,\,\, \forall t\neq 0,
  \]
  where
  \[
  \Phi(t) = \int_{0}^{|t|} \phi(s)sds.
  \]
\item\label{H3} The function $\displaystyle\frac{\phi(t)}{t^{m-2}}$ is nonincreasing in $(0, +\infty)$.

\item\label{H4} The function $\phi$ is monotone.

\item\label{H5} There exists a constant $c>0$ such that
\[
\vert \phi^{'}(t)t\vert \leq c\phi(t), \quad \forall \ t \in [0, +\infty).
\]

\end{enumerate}

Hereafter, we will say that $\Phi \in \mathcal{C}_m$ if
$$
\Phi(t) \geq |t|^{m}, \,\,\,\,\,\, \forall t \in \mathbb{R}. \leqno{(\mathcal{C}_m)}
$$
Moreover, let us denote by $\gamma$ the following real number
$$
\gamma=
\left\{
\begin{array}{l}
m, \,\,\,\, \mbox{if} \,\,\, \Phi \in \mathcal{C}_m, \\
\mbox{}\\
l, \,\,\,\, \mbox{if} \,\,\, \Phi \notin \mathcal{C}_m.
\end{array}
\right.
$$

Here, we would like to detach that the functions $\phi$ associated with each N-function mentioned in this introduction, fulfill the conditions $(\phi_1)$-$(\phi_5)$.

\begin{flushleft}
\textbf{Conditions on $f$:}
\end{flushleft}
The function $f: \mathbb{R}\rightarrow \mathbb{R}$ is a $C^{1}$- function verifying:
\begin{enumerate}[label={($f_\arabic{*}$})]
  \item\label{f1}  There are functions $r,b:[0,+\infty) \to [0,+\infty)$ such that
  \[
  \limsup_{|t|\rightarrow 0} \frac{f^{'}(t)}{(r(|t|)|t|)^{'}}=0 \quad \mbox{and} \quad \limsup_{|t|\rightarrow +\infty} \frac{|f^{'}(t)|}{(b(|t|)|t|)^{'}}<+\infty.
  \]
  \item\label{f2} There exists $\theta > m$ such that
  \[
  0< \theta F(t) \leq f(t)t, \quad \forall t>0
  \]
  where
  \[
  F(t)= \int_{0}^{s}f(s)ds.
  \]
  \item\label{f3}  The function $\displaystyle\frac{f(t)}{t^{m -1}}$ is increasing for $t>0$.
\end{enumerate}
Related to functions $r $ and $b$, we assume that they are $C^{1}$- functions  satisfying the following conditions:
\begin{enumerate}[label={($r_\arabic{*}$})]
\setcounter{enumi}{0}
\item\label{R1} $r$ is increasing.
\item\label{R2} There exists a constant $\overline{c}>0$ such that
\[
\vert r^{'}(t)t \vert \leq \overline{c} r(t), \quad \forall \ t\geq 0.
\]
\item\label{R3} There exist positive constants $r_{1}$ and $r_{2}$ such that
 \[
  r_{1}\leq \frac{r(t)t^{2}}{R(t)}\leq r_{2}, \quad \forall t>0,
  \]
  where
  \[
  R(t) = \int_{0}^{|t|} r(s)sds.
  \]
\item\label{R4} The function $R$ satisfies
\[
  \limsup_{t\rightarrow 0}\frac{R(t)}{\Phi(t)}< +\infty \quad \mbox{and} \quad \limsup_{|t|\rightarrow +\infty}\frac{R(t)}{\Phi_{*}(t)}=0.
  \]

\end{enumerate}

\begin{enumerate}[label={($b_\arabic{*}$})]
\setcounter{enumi}{0}
\item\label{B1} $b$ is increasing.
 \item\label{B2} There exists a constant $\widetilde{c}>0$ such that
\[
\vert b^{'}(t)t\vert \leq \widetilde{c}b(t), \quad \forall t\geq 0.
\]
\item\label{B3} There exist positive constants $b_{1}, b_{2} \in (1, \gamma^*)$ verifying
 \[
  b_{1}\leq \frac{b(t)t^{2}}{B(t)}\leq b_{2}, \,\,\, \forall t> 0,
  \]
  where $\gamma^*=\frac{N\gamma}{N-\gamma}$ and
  \[
  B(t) = \int_{0}^{|t|} b(s)sds.
  \]
\item\label{B4} The function $B$ satisfies
\[
  \limsup_{t\rightarrow 0}\frac{B(t)}{\Phi(t)}< +\infty \quad \mbox{and} \quad \limsup_{|t|\rightarrow +\infty}\frac{B(t)}{\Phi_{*}(t)}=0,
  \]
\end{enumerate}
where $\Phi_{*}$ is the Sobolev conjugate function, which is defined by inverse function of
  \[
  G_{\Phi}(t) = \int_{0}^{t}\frac{\Phi^{-1}(s)}{s^{1+\frac{1}{N}}}ds.
  \]

Using the above hypotheses, we are able to state our main result.

\begin{thm}\label{T1}
Suppose that  \ref{H1}-\ref{H5}, \ref{R1}-\ref{R4}, \ref{B1}-\ref{B4}, \ref{f1}-\ref{f3}, $(V_0)$-$(V_1)$ hold. Then, for any $\delta > 0$ small enough, there exists  $\epsilon_{\delta}> 0$ such that $(P_{\epsilon} )$ has at least $cat_{M_{\delta}}(M)$ positive solutions, for any $0 < \epsilon < \epsilon_{\delta}$, where
$$
M = \big\{x \in \Omega \ : \ V(x)=V_0 \big\}
$$
and
\[
M_{\delta} = \big\{x \in \mathbb{R}^{N} \ : \ dist(x,M) \leq \delta \big\}.
\]
Moreover, if $u_{\epsilon}$ denotes one of these solutions and $x_{\epsilon} \in \mathbb{R}^{N}$ is its global maximum, we have that
\[
\lim_{\epsilon \rightarrow 0}V(\epsilon x_{\epsilon}) = V_{0}.
\]
\end{thm}

We would like point out that, if $Y$ is a closed subset of a topological space $X$, the Lusternik-Schnirelman category $cat_{X}(Y)$ is the least number of closed and contractible sets in $X$ which cover $Y$.

The plan of the paper is as follows: In Section 2, we will prove the existence and multiplicity of solution for an auxiliary problem, more precisely, using Lusternik-Schnirelman category theory, we show that the auxiliary problem has at last $cat_{M_{\delta}}(M)$ positive solutions for $\epsilon$ small enough. In Section 3, we make some estimates to prove that the solutions found  for the auxiliary problem in Section 2 are solutions for the original problem. Finally, we write  an Appendix A, where we show the existence of a special function used in Section 2.

\section{An auxiliary problem}
In this is section, motivated for some arguments explored in Alves and Figueiredo \cite{AF2005}, and mainly in del Pino and Felmer \cite{DF1}, we will show the existence and  multiplicity of positive solutions for an auxiliary problem. To this end, we need to fix some notations, however if the reader does not know the main properties involving the Orlicz-Sobolev spaces, we suggest to read the Section 2 in \cite{AS} for a brief review, and for a more complete study, see \cite{Adams1}, \cite{Adams2}, \cite{DT}, \cite{MU} and \cite{Rao}.

Since we intend to find positive solutions, we will assume that
\begin{eqnarray}\label{posf}
f(t) = 0 \quad \mbox{for all} \ t< 0.
\end{eqnarray}

Let $\theta$ be the number given in \ref{f3}, $a, k > 0$ satisfying
$$
k>\displaystyle\frac{(\theta - l)}{(\theta - m)}\frac{m}{l} \quad \mbox{and} \quad \displaystyle\frac{f(a)}{\phi(a)a} = \frac{V_{0}}{k}.
$$
Using the above numbers, let us define the function
\begin{align}
\widehat{f}(s) = \left\{
\begin{array}
[c]{rcl}%
f(s) && if ~  s\leq a\\
\displaystyle\frac{V_{0}}{k}\phi(s)s  && if ~ s>a.
\end{array}\nonumber
\right.%
\end{align}
Fixing $t_0 < a < t_1$ with $t_0,t_1 \approx a$, it is possible to build a function $\eta \in C^1([t_0, t_1])$ satisfying
\begin{enumerate}[label={($\eta_\arabic{*}$})]
\setcounter{enumi}{0}
\item\label{eta1} $\eta(s) \leq \widehat{f}(s)$ for all $s \in [t_0, t_1]$,
  \item\label{eta2} $\eta(t_0) = \widehat{f}(t_0)$ and $\eta(t_1)= \widehat{f}(t_1)$,
  \item\label{eta3} $\eta^{'}(t_0) = (\widehat{f})^{'}(t_0)$ and $\eta^{'}(t_1)= (\widehat{f})^{'}(t_1)$,
  \item\label{eta4} The function $s \rightarrow \displaystyle\frac{\eta(s)}{\phi(s)s}$ is nondecreasing for all $s \in [t_0, t_1]$.
\end{enumerate}

The reader can find the details involving the existence of  $\eta$ in  Appendix  A.

Using the functions $\eta$ and $\widehat{f}$, let us consider two new functions
\begin{align}
\widetilde{f}(s) = \left\{
\begin{array}
[c]{rcl}%
\widehat{f}(s) && if ~  s \not\in [t_0, t_1],\\
\eta(s)  && if ~ s \in [t_0, t_1],
\end{array}\nonumber
\right.%
\end{align}
and
\[
g(x, s) = \chi_{\Omega}(x)f(s) + (1-\chi_{\Omega}(x))\widetilde{f}(s),
\]
where $\chi_{\Omega}$ is the characteristic function related to the set $\Omega$. From definition of $g$, we see that  $g$ is a Carath\'eodory function verifying
\begin{eqnarray}\label{posg}
g(x, s) = 0, \quad \forall (x,s) \in \mathbb{R}^{N} \times (-\infty, 0]
\end{eqnarray}
and
\begin{eqnarray}\label{gf}
g(x, s) \leq f(s), \quad \forall (x,s) \in \mathbb{R}^{N} \times \mathbb{R}.
\end{eqnarray}
Moreover, for each $x \in \mathbb{R}^{N}$, the function $s\rightarrow g(x, s)$ is of class $C^{1}$ and it satisfies the following conditions:
\begin{enumerate}[label={($g_\arabic{*}$})]
	\setcounter{enumi}{0}
	\item\label{g1} $\displaystyle\limsup_{|s|\rightarrow 0}\frac{g(x, s)}{\phi(|s|)|s|}=0$, uniformly in $x \in \mathbb{R}^{N}$.

     \item\label{g2} $\displaystyle\limsup_{|s|\rightarrow +\infty}\frac{ g(x, s)}{b(|s|)|s|}< +\infty$, uniformly in $x \in \mathbb{R}^{N}$.
	
	\item\label{g3}$0\leq\theta G(x, s) = \theta\displaystyle\int_{0}^{s}g(x, t)dt \leq g(x, s)s$, $\forall (x,s) \in \Omega \times (0,+\infty)$.
	\item\label{g4} $0<lG(x, s)\leq g(x, s)s\leq \displaystyle\frac{V_0}{k}\phi(s)s^{2}$, $\forall (x,s) \in \Omega^{c} \times (0,+\infty)$.
	\item\label{g5}The function $s \rightarrow \displaystyle\frac{g(x, s)}{\phi(s)s}$ is nondecreasing for each $x \in \mathbb{R}^{N}$ and for all $s > 0$.
\end{enumerate}

Using the function $g$, we can consider the auxiliary problem
\begin{align}
\left\{
\begin{array}
[c]{rcl}%
- \Delta_{\Phi}u + V(\epsilon x)\phi(|u|)u & = & g(\epsilon x, u)~  \mbox{in}~ \mathbb{R}^{N},\\
u \in W^{1, \Phi}(\mathbb{R}^{N}). &  &
\end{array}
\right. \tag{$ \widetilde{P}_{\epsilon} $}\label{PA2}%
\end{align}

\vspace{0.5 cm}

Here, we would like to point out that if $\Omega_\epsilon$ denotes the set
$$
\Omega_\epsilon=\{x \in \mathbb{R}^{N}\,:\,\epsilon x \in \Omega\}=\Omega / \epsilon
$$
and $u$ is a positive solution of \eqref{PA2} with $u(x)\leq t_0$ for all $x \in \mathbb{R}^{N}\backslash \Omega_\epsilon$, then $u$ is also a positive solution of \eqref{P2}.

\subsection{Preliminary results}


In what follows, let us denote by $J_{\epsilon} : X_{\epsilon} \rightarrow \mathbb{R}$ the energy functional related to \eqref{PA2} given by
\[
J_{\epsilon}(u)= \int_{\mathbb{R}^{N}}\Phi(\vert\nabla u\vert)dx + \int_{\mathbb{R}^{N}}V(\epsilon x)\Phi(\vert u\vert)dx - \int_{\mathbb{R}^{N}}G(\epsilon x, u)dx,
\]
where $X_{\epsilon}$ denotes the subspace of $W^{1, \Phi}(\mathbb{R}^{N})$ given by
\[
X_{\epsilon}= \Big\{ u \in W^{1, \Phi}(\mathbb{R}^{N})\, :\, \int_{\mathbb{R}^{N}}V(\epsilon x)\Phi(\vert u \vert)dx < +\infty\Big\},
\]
endowed with the norm
\[
\Vert u\Vert_{\epsilon} = \Vert \nabla u\Vert_{\Phi} + \Vert u\Vert_{\Phi, V_{\epsilon}},
\]
where
$$
\Vert \nabla u \Vert_{\Phi}:= \inf\Big\{ \lambda > 0;  \int_{\mathbb{R}^{N}}\Phi\Big(\frac{\vert \nabla u\vert}{\lambda}\Big)dx \leq 1\Big\}
$$
and
\[
\Vert u\Vert_{\Phi, V_{\epsilon}}:= \inf\Big\{ \lambda > 0;  \int_{\mathbb{R}^{N}}V(\epsilon x)\Phi\Big(\frac{\vert u\vert}{\lambda}\Big)dx \leq 1\Big\}.
\]
From $(V_{0})$, it follows that the embeddings
\[
X_{\epsilon}\hookrightarrow L^{\Phi}( \mathbb{R}^{N}) \quad \mbox{and}\quad X_{\epsilon}\hookrightarrow L^{B}( \mathbb{R}^{N})
\]
are continuous. Using the above embeddings, a direct computation yields   $J_{\epsilon} \in C^{1}(X_{\epsilon}, \mathbb{R})$ with
\[
J^{'}_{\epsilon}(u)v = \int_{\mathbb{R}^{N}}\phi(\vert\nabla u\vert)\nabla u \nabla v \, dx + \int_{\mathbb{R}^{N}}V(\epsilon x)\phi(\vert u\vert)uv \, dx - \int_{\mathbb{R}^{N}}g(\epsilon x, u)v \, dx,
\]
for all $u,v \in X_{\epsilon}$. Thereby, $u \in X_\epsilon$ is a weak solution of \eqref{PA2} if, and only if,
$u$ is a critical point of $J_{\epsilon}$. Furthermore, by \eqref{posg}, the critical points of $J_{\epsilon}$ are
nonnegative.
\begin{lem}\label{Lmtd21}
Let $(u_{n})$ be a sequence $(PS)_{c}$. Then, $(u_{n})$ is a bounded sequence in $X_{\epsilon}$.
\end{lem}
\begin{pf}
Since $(u_{n})$ is a $(PS)_{c}$ sequence for $J_{\epsilon}$, there is $C_{1}>0$ such that
\begin{eqnarray*}
J_{\epsilon}(u_{n}) - \frac{1}{\theta}J_{\epsilon}^{'}(u_{n})u_{n}\leq C_{1}(1 + \Vert u_{n}\Vert_{\epsilon}), \quad \forall n \in \mathbb{N}.
\end{eqnarray*}
On the other hand, by $(\phi_2)$ and $(g_3)-(g_4)$,
$$
J_{\epsilon}(u_{n}) - \frac{1}{\theta}J_{\epsilon}^{'}(u_{n})u_{n} \geq C\Bigg(\int_{\mathbb{R}^{N}}\Phi(|\nabla u_{n}|)dx + \int_{\mathbb{R}^{N}}V(\epsilon x)\Phi(u_{n})dx\Bigg),
$$
where $C = \Bigg[\Big(1-\frac{m}{\theta}\Big)-\Big(1-\frac{l}{\theta}\Big)\frac{m}{kl}\Bigg]>0$. Hence, by \cite[Lemma 2.1]{FN},
\begin{eqnarray*}
J_{\epsilon}(u_{n}) - \frac{1}{\theta}J_{\epsilon}^{'}(u_{n})u_{n} &\geq & C_{1}\big( \xi_{0}(\Vert \nabla u_{n}\Vert_{\Phi}) + \xi_{0}(\Vert u_{n}\Vert_{\Phi, V_{\epsilon}})\big), \quad \forall n \in \mathbb{N}.
\end{eqnarray*}
Now, the proof follows as in \cite[Lemma 4.2]{AGJ}.
\end{pf}
\begin{lem}\label{Lema22}
Let $(u_n)$ be a $(PS)_d$ sequence for $J_\epsilon$. Then for each $\tau$, there exists
$\rho_{0} = \rho_{0}(\tau) > 0$ such that
\[
\limsup_{n \rightarrow +\infty}\int_{\mathbb{R}^{N}\backslash B_{\rho_0}(0)}\Big[ \Phi(\vert u_n\vert) + V(\epsilon x)\Phi(\vert u_n\vert)\Big]dx < \tau.
\]
\end{lem}

\begin{pf}
For each $\rho>0$, let $\xi_\rho \in C^{\infty}(\mathbb{R}^{N})$ verifying
\[
\xi_{\rho}(x) =
\left\{
\begin{array}{rcl}
  0, \quad x \in B_{\frac{\rho}{2}}(0)\\
  1,  \quad  x \notin B_{\rho}(0)\\
\end{array}
\right.
\]
with $0\leq \xi_{\rho}(x)\leq 1$ and $\vert \nabla \xi_{\rho}\vert \leq \displaystyle\frac{C}{\rho}$, where $C$ is a constant independent of $\rho$. Note that
\begin{eqnarray*}
J^{'}_{\epsilon}(u_n)(\xi_\rho u_n) &=& \int_{\mathbb{R}^{N}}\phi(\vert \nabla u_n\vert)\nabla u_{n}\nabla(\xi_\rho u_n)dx + \int_{\mathbb{R}^{N}}V(\epsilon x)\phi(\vert u_n\vert)u_{n}^{2}\xi_\rho dx\\
 && - \int_{\mathbb{R}^{N}}g(\epsilon x, u_n)u_n\xi_\rho dx.
\end{eqnarray*}
Choosing $\rho>0$ such that $\Omega_{\epsilon} \subset B_{\frac{\rho}{2}}(0)$, the condition \ref{H2} ensures that
\begin{eqnarray*}
l\int_{\mathbb{R}^{N}}\xi_\rho\Big[\Phi(\vert \nabla u_n\vert) + V(\epsilon x)\Phi(\vert u_n\vert)\Big]dx &\leq& J^{'}_{\epsilon}(u_n)(\xi_\rho u_n)- \int_{\mathbb{R}^{N}}u_{n}\phi(\vert \nabla u_n\vert)\nabla u_{n}\nabla \xi_\rho dx\\
&&+ \int_{\mathbb{R}^{N}\backslash \Omega_{\epsilon}}g(\epsilon x, u_n)u_n\xi_\rho dx.
\end{eqnarray*}
Gathering \ref{g4} and \ref{H2},
\begin{eqnarray*}
l\int_{\mathbb{R}^{N}}\xi_\rho\Big[\Phi(\vert \nabla u_n\vert) + V(\epsilon x)\Phi(\vert u_n\vert)\Big]dx &\leq& J^{'}_{\epsilon}(u_n)(\xi_\rho u_n)- \int_{\mathbb{R}^{N}}u_{n}\phi(\vert \nabla u_n\vert)\nabla u_{n}\nabla \xi_\rho dx\\
&&+ \frac{m}{k}\int_{\mathbb{R}^{N}}V(\epsilon x)\Phi(\vert u_{n}\vert)\xi_\rho dx.
\end{eqnarray*}
Since $(\xi_\rho u_n)$ is bounded in $X_{\epsilon}$ and $k>\displaystyle\frac{m}{l}$, by H\"{o}lder inequality there exists a constant $C_1>0$ such that
\begin{eqnarray*}
\int_{\mathbb{R}^{N}}\xi_\rho\Big[\Phi(\vert \nabla u_n\vert) + V(\epsilon x)\Phi(\vert u_n\vert)\Big]dx &\leq& o_{n}(1)+ \frac{C_1}{\rho}.
\end{eqnarray*}
Now, fixing $\tau>0$, there exists $\rho_0>0$ such that $\frac{C_{1}}{\rho_0}<\tau$. Then,
\begin{eqnarray*}
\int_{\mathbb{R}^{N}\backslash B_{\rho_0}(0)}\Big[\Phi(\vert \nabla u_n\vert) + V(\epsilon x)\Phi(\vert u_n\vert)\Big]dx &\leq& o_{n}(1)+ \tau.
\end{eqnarray*}
Passing to the limit in the last inequality, it follows that
\[
\limsup_{n \rightarrow +\infty}\int_{\mathbb{R}^{N}\backslash B_{\rho_0}(0)}\Big[\Phi(\vert u_n\vert) + V(\epsilon x)\Phi(\vert u_n\vert)\Big]dx < \tau.
\]
\end{pf}
The Lemma below establishes an important property involving the $(PS)$ sequences for $J_{\epsilon}$. Since it follows repeating the same arguments used in \cite[Lemma 4.3]{AGJ}, we will omit its proof.
\begin{lem}\label{LemaP}
Let $(u_n)$ be a $(PS)_d$ sequence for $J_\epsilon$ with $u_{n}\rightharpoonup u$ in $X_{\epsilon}$. Then,
\begin{eqnarray}\label{cgrad}
\nabla u_{n}(x)\rightarrow \nabla u(x)\quad \mbox{a. e.  in}\quad \mathbb{R}^{N}.
\end{eqnarray}
As a consequence of the above limit, we deduce that $u$ is a critical point for $J_\epsilon$.
\end{lem}
\begin{prop}\label{prop2.1}
The functional $J_\epsilon$ verifies the $(PS)$ condition.
\end{prop}
\begin{pf}
Let $(u_n)$ be a $(PS)_c$ sequence for $J_{\epsilon}$. From Lemma \ref{Lmtd21}, there exists $u \in X_{\epsilon}$ such that
\begin{eqnarray}\label{In}
u_n \rightharpoonup u \, \,\, \mbox{in} \,\,  X_{\epsilon}.
\end{eqnarray}
First of all, we know by Lemma \ref{Lema22} that fixed $\tau>0$, there is $\rho_0>0$ such that
\[
\limsup_{n \rightarrow +\infty}\int_{\mathbb{R}^{N}\backslash B_{\rho_0}(0)}\Big[\Phi(\vert u_n\vert) + V(\epsilon x)\Phi(\vert u_n\vert)\Big]dx < \tau.
\]
Increasing $\rho_0>0$ if necessary, the above limit together with $\Delta_2$-condition gives
\begin{equation} \label{gradiente}
\limsup_{n \rightarrow +\infty}\int_{\mathbb{R}^{N}}\Phi(\vert \nabla u_n - \nabla u\vert)dx \leq \limsup_{n \rightarrow +\infty}\int_{B_{\rho_0}(0)}\Phi(\vert \nabla u_n - \nabla u\vert)dx + 2\tau\quad
\end{equation}
and
\begin{eqnarray*}
\limsup_{n \rightarrow +\infty}\int_{\mathbb{R}^{N}}V(\epsilon x)\Phi(\vert  u_n - u\vert)dx &\leq& \limsup_{n \rightarrow +\infty}\int_{B_{\rho_0}(0)}V(\epsilon x)\Phi(\vert u_n -  u\vert)dx + 2\tau.
\end{eqnarray*}
By \eqref{In}, up to a subsequence, $u_n \rightarrow u$ in $L^{\Phi}(B_{\rho_0}(0))$. This information combined with the last limit guarantees that
$$
\limsup_{n \rightarrow +\infty}\int_{\mathbb{R}^{N}}V(\epsilon x)\Phi(\vert  u_n - u\vert)dx \leq 2\tau.
$$
As $\tau$ is arbitrary, we conclude that
\begin{eqnarray}\label{cg4}
\limsup_{n \rightarrow +\infty}\int_{\mathbb{R}^{N}}V(\epsilon x)\Phi(\vert  u_n - u\vert)dx=0.
\end{eqnarray}

Now, we will show that
\begin{eqnarray*}
\limsup_{n \rightarrow +\infty}\int_{B_{\rho_0}(0)}\Phi(\vert \nabla u_n - \nabla u\vert)dx =0.
\end{eqnarray*}
By Lemma \ref{LemaP},
\[
\Phi(\vert \nabla u_n(x) -\nabla u(x)\vert) \rightarrow 0 \quad \mbox{a. e. in }\, B_{\rho_0}(0).
\]
Moreover, from $\Delta_2$-condition and \ref{H2}, there exist constants $c_1, c_2>0$ such that
\[
\Phi(\vert \nabla u_n -\nabla u\vert) \leq c_1\phi(\vert\nabla u_n \vert)\vert\nabla u_n \vert^{2} + c_2\phi(\vert\nabla u \vert)\vert\nabla u \vert^{2}.
\]
Using again Lemma \ref{LemaP}
\[c_1\phi(\vert\nabla u_n \vert)\vert\nabla u_n \vert^{2} + c_2\phi(\vert\nabla u \vert)\vert\nabla u \vert^{2} \rightarrow (c_1+c_2)\phi(\vert\nabla u \vert)\vert\nabla u \vert^{2} \quad \mbox{a. e. in }\, B_{\rho_0}(0).
\]
On the other hand, as in \cite[Lemma 4.3]{AGJ}, 
\[
\int_{B_{\rho_{0}}(0)}(\phi(\vert \nabla u_n\vert)\nabla u_n - \phi(\vert\nabla u\vert)\nabla u)(\nabla u_{n} - \nabla u)dx = o_{n}(1),
\]
and so,
\[
\int_{B_{\rho_0}(0)}\phi(\vert \nabla u_n\vert)|\nabla u_n|^{2}  dx\rightarrow\int_{B_{\rho_0}(0)}\phi(\vert \nabla u \vert)\vert\nabla u\vert^{2} dx.
\]
Therefore,
\[
\int_{B_{\rho_0}(0)}\Big[c_1\phi(\vert \nabla u_n\vert)\vert \nabla u_n\vert^{2} + c_2\phi(\vert \nabla u\vert)\vert\nabla u\vert^{2}\Big]dx \rightarrow (c_1 + c_2)\int_{B_{\rho_0}(0)}\phi(\vert \nabla u\vert)\vert\nabla u\vert^{2}dx.
\]
Applying the general Lebesgue's theorem, we infer that
\[
\lim_{n\rightarrow +\infty}\int_{B_{\rho_0}(0)}\Phi(\vert \nabla u_n -\nabla u\vert)dx = 0.
\]
Combining the last limit with \eqref{gradiente}, we obtain
\begin{eqnarray}\label{cg5}
\lim_{n\rightarrow +\infty}\int_{\mathbb{R}^{N}}\Phi(\vert \nabla u_n -\nabla u\vert)dx = 0.
\end{eqnarray}
From \eqref{cg4} and \eqref{cg5},
\[
u_n\rightarrow u\quad \mbox{in}\, \, X_{\epsilon},
\]
showing that $J_{\epsilon}$ verifies the $(PS)$ condition.
\end{pf}

Next, we will show some results involving $J_{\epsilon}$ and its Nehari manifold.  We recall that the Nehari manifold associated to $J_{\epsilon}$ is given by,
\[
\mathcal{N}_{\epsilon} = \big\{ u \in X_{\epsilon}\backslash \{0 \} \ : \ J^{'}_{\epsilon}(u)u = 0 \big\}.
\]

Our first lemma shows that the Nehari manifold has a positive distance from the origin in $X_\epsilon$. Since it follows by using standard arguments, we omit its proof.
\begin{lem} \label{Nehari2}
	For all $u \in \mathcal{N}_{\epsilon}$, there exists $\sigma>0$, which is independent of $\epsilon$, such that
	\begin{eqnarray*}
		\Vert u\Vert_{\epsilon} > \sigma.
	\end{eqnarray*}
\end{lem}

The next lemma establishes an important characterization involving the mountain pass level, which is useful in a lot of problems. In what follows,  we denote by $c_{\epsilon, 1}$ and $c_{\epsilon, 2}$ the following numbers
\[
c_{\epsilon, 1} = \inf_{u \in \mathcal{N}_{\epsilon}}J_{\epsilon}(u) \quad \mbox{and} \quad c_{\epsilon, 2} = \inf_{u \in X_{\epsilon}\backslash\{0\}}\max_{t\geq 0}J_{\epsilon}(tu).
\]
Fixing the  subset
\[
\mathcal{A}_{\epsilon} = \big\{ u \in X_{\epsilon} \, \, :\, \, u^{+}\neq0 \quad \mbox{and}\quad |supp(u)\cap\Omega_{\epsilon}|>0\big\}
\]
and the number
\[
\overline{c}_{\epsilon, 2} = \inf_{u \in \mathcal{A}_{\epsilon}}\max_{t\geq 0}J_{\epsilon}(tu),
\]
it is easy to see that
\[
c_{\epsilon, 2} = \overline{c}_{\epsilon, 2}.
\]

\begin{lem}\label{Lcar}
	Assume that \ref{H1}-\ref{H3}, \ref{R1}-\ref{R4}, \ref{B1}-\ref{B4}, \ref{f1}-\ref{f3}  and $(V_0)$-$(V_1)$ hold. Then, for each $u \in \mathcal{A}_{\epsilon}$, there exists a unique $t_{u}>0$ such that $t_{u}u \in \mathcal{N}_{\epsilon}$ and $J_{\epsilon}(t_{u}u) = \displaystyle\max_{t \geq 0}J_{\epsilon}(tu)$. Moreover,
	\[
	c_{\epsilon} = c_{\epsilon, 1} = c_{\epsilon, 2},
	\]
	where $c_{\epsilon}$ denotes the mountain pass level associated with $J_{\epsilon}$.
\end{lem}
\begin{pf}
	For each  $u \in \mathcal{A}_{\epsilon}$, we define $h_{\epsilon}(t)= J_{\epsilon}(tu)$, that is,
	\begin{eqnarray*}
		h_{\epsilon}(t) = \displaystyle\int_{\mathbb{R}^{N}}\Phi(\vert\nabla (tu)\vert)dx + \displaystyle\int_{\mathbb{R}^{N}}V(\epsilon x)\Phi(\vert tu\vert)dx - \displaystyle\int_{\mathbb{R}^{N}}G(\epsilon x, tu)dx.
	\end{eqnarray*}
	
	\begin{flushleft}
		\textbf{Existence}
	\end{flushleft}
	By a direct computation,  $h_{\epsilon}(t)> 0$ for $t$ enough small and $h_{\epsilon}(t)<0$ for $t$ sufficiently large. Thus, there is $t_{u}>0$ such that
	\[
	h_{\epsilon}(t_{u}) = \max_{t\geq 0}h_{\epsilon}(t) = \max_{t\geq 0}J_{\epsilon}(tu),
	\]
	implying that $h_{\epsilon}^{'}(t_{u}) = 0$, that is, $J^{'}_{\epsilon}(t_{u}u)u = 0$, and so, $t_{u}u \in \mathcal{N}_{\epsilon}$.
	\begin{flushleft}
		\textbf{Uniqueness}
	\end{flushleft}
	Suppose that there exist $t_{1}, t_{2} >0$ such that $t_{1}u, t_{2}u \in \mathcal{N}_{\epsilon}$ and $t_{1}<t_{2}$. Then,
	\begin{eqnarray*}
		\int_{\mathbb{R}^{N}}\phi(|\nabla (t_{1}u)|)|\nabla(t_{1}u)|^{2}dx &+&  \int_{\mathbb{R}^{N}}V(\epsilon x)\phi(|t_{1}u|)|t_{1}u|^{2}dx = \int_{[u>0]}g(\epsilon x, t_{1}u)t_{1}udx\\
	\end{eqnarray*}
	and
	\begin{eqnarray*}
		\int_{\mathbb{R}^{N}}\phi(|\nabla (t_{2}u)|)|\nabla(t_{2}u)|^{2}dx &+&  \int_{\mathbb{R}^{N}}V(\epsilon x)\phi(|t_{2}u|)|t_{2}u|^{2}dx = \int_{[u>0]}g(\epsilon x, t_{2}u)t_{2}udx.\\
	\end{eqnarray*}
	Setting $\upsilon(t) = \displaystyle\frac{\phi(t)}{t^{m-2}}$ for all $t>0$, we have
	\begin{eqnarray*}
		&&\displaystyle\int_{\mathbb{R}^{N}}\big(\upsilon(\vert t_1|\nabla u|\vert) - \upsilon(\vert t_2|\nabla u|\vert)\big)|\nabla u|^{m}dx + \displaystyle\int_{\mathbb{R}^{N}}V(\epsilon x)\big(\upsilon(\vert t_1| u|\vert) - \upsilon(\vert t_2|u|\vert)\big)|u|^{m}dx\\
		&=& \displaystyle\int_{[u>0]}\Bigg[\displaystyle\frac{g(\epsilon x, t_{1}u)}{(t_{1}u)^{m-1}}-\displaystyle\frac{g(\epsilon x, t_{2}u)}{(t_{2}u)^{m-1}}\Bigg]u^{m}dx.
	\end{eqnarray*}
	Thus, from $(V_{0})$ and \ref{H3},
	\begin{eqnarray*}
		&&\displaystyle\int_{\mathbb{R}^{N}}\big(\upsilon(\vert t_1|\nabla u|\vert) - \upsilon(\vert t_2|\nabla u|\vert)\big)|\nabla u|^{m}dx\\
		&&+ \frac{V_{0}}{k}\displaystyle\int_{(\mathbb{R}^{N}\backslash\Omega_{\epsilon})\cap[u>0]}\big(\upsilon(\vert t_1| u|\vert) - \upsilon(\vert t_2|u|\vert)\big)|u|^{m}dx\\
		&&\leq\displaystyle\int_{\Omega_{\epsilon}\cap[u>0]}\Bigg[\displaystyle\frac{f(t_{1}u)}{(t_{1}u)^{m-1}}-\displaystyle\frac{f(t_{2}u)}{(t_{2}u)^{m-1}}\Bigg]u^{m}dx \\
		&&+\displaystyle\int_{(\mathbb{R}^{N}\backslash\Omega_{\epsilon})\cap[u>0]}\Bigg[\displaystyle\frac{\widetilde{f}(t_{1}u)}{(t_{1}u)^{m-1}}-\displaystyle\frac{\widetilde{f}(t_{2}u)}{(t_{2}u)^{m-1}}\Bigg]u^{m}dx,
	\end{eqnarray*}
	implying that
	\begin{eqnarray*}
		&&\displaystyle\int_{\mathbb{R}^{N}}\big(\upsilon(\vert t_1|\nabla u|\vert) - \upsilon(\vert t_2|\nabla u|\vert)\big)|\nabla u|^{m}dx\\
		&&+ \displaystyle\int_{(\mathbb{R}^{N}\backslash\Omega_{\epsilon})\cap[u>0]}\Bigg[\bigg(\displaystyle\frac{V_{0}}{k}\upsilon(\vert t_1| u|\vert) -\displaystyle\frac{\widetilde{f}(t_{1}u)}{(t_{1}u)^{m-1}}\bigg) - \bigg(\displaystyle\frac{V_{0}}{k}\upsilon(\vert t_2| u|\vert)-\displaystyle\frac{\widetilde{f}(t_{2}u)}{(t_{2}u)^{m-1}}\bigg)\Bigg]u^{m}dx\\
		&\leq&\displaystyle\int_{\Omega_{\epsilon}\cap[u>0]}\Bigg[\displaystyle\frac{f(t_{1}u)}{(t_{1}u)^{m-1}}-\displaystyle\frac{f(t_{2}u)}{(t_{2}u)^{m-1}}\Bigg]u^{m}\,dx.
	\end{eqnarray*}
	Now, consider the function
	\[
	h(t) = \displaystyle\frac{V_{0}}{k}\upsilon(t) - \frac{\widetilde{f}(t)}{t^{m-1}},
	\]
	and note that,
	\[
	h(t) = \upsilon(t)h_{1}(t),
	\]
	where
	$$
	h	_{1}(t) = \displaystyle\frac{V_{0}}{k}-\displaystyle\frac{\widetilde{f}(t)}{\phi(t)t}.
	$$
	
	As $\upsilon$, $h_{1}$ are non increasing and nonnegative, $h$ is non increasing. Hence,  $h(t_{1}u)\geq h(t_{2}u)$, and so,
	\begin{eqnarray*}
		0&\leq&\displaystyle\int_{\mathbb{R}^{N}}\big(\upsilon(\vert t_1|\nabla u|\vert) - \upsilon(\vert t_2|\nabla u|\vert)\big)|\nabla u|^{m}dx +\displaystyle\int_{(\mathbb{R}^{N}\backslash\Omega_{\epsilon})\cap[u>0]}\big(h(t_{1}u) - h(t_{2}u)\big)u^{m}dx\\
		&\leq&\displaystyle\int_{\Omega_{\epsilon}\cap[u>0]}\Bigg[\displaystyle\frac{f(t_{1}u)}{(t_{1}u)^{m-1}}-\displaystyle\frac{f(t_{2}u)}{(t_{2}u)^{m-1}}\Bigg]u^{m}dx<0,
	\end{eqnarray*}
	which is an absurd. Therefore, $t_{1} = t_{2}$. Now, the proof follows arguing as in \cite[Theorem 4.2]{Willen}.
	\end{pf}

\subsection{Study of the $(PS)$ condition for $J_\epsilon$ on $\mathcal{N}_{\epsilon}$}

In this subsection, our main goal is to study the $(PS)$ condition for $J_\epsilon$ on $\mathcal{N}_{\epsilon}$, which will be used later on. In what follows, without loss of generality, we will assume that
$$
V(0)=\min_{z \in \mathbb{R}^{N}}V(z)=V_0.
$$
Using the above number, we fix the following quasilinear problem
\begin{align}
\left\{
\begin{array}{l}
- \Delta_{\Phi}u + V_0\phi(\vert u\vert)u = f(u)~  \mbox{in}~ \mathbb{R}^{N},\\
u \in W^{1, \Phi}(\mathbb{R}^{N}).
\end{array}
\right. \tag{$ P_{0} $}\label{PAu2}%
\end{align}
We recall that the weak solutions of \eqref{PAu2} are critical points of the functional
\[
E_{0}(u) = \int_{\mathbb{R}^{N}}\Phi(\vert \nabla u \vert)dx + V_0\int_{\mathbb{R}^{N}}\Phi(\vert u \vert)dx -\int_{\mathbb{R}^{N}}F(u)dx,
\]
which is well defined in $Y= W^{1, \Phi}(\mathbb{R}^{N})$ endowed with the norm
\[
\Vert u\Vert_{Y} = \Vert \nabla u\Vert_{\Phi} + V_0 \Vert u\Vert_{\Phi}.
\]
Moreover, let us denote by $d_{0}$ the mountain pass level  of $E_{0}$, and by $\mathcal{M}_{0}$, the Nehari manifold
\[
\mathcal{M}_{0} = \big\{u \in Y \backslash \{0\}\,:\, E^{'}_{0}(u)u=0 \big\}.
\]

The next lemma will be used to prove that  $J_\epsilon$ verifies the $(PS)$ condition on $\mathcal{N}_{\epsilon}$.

\begin{lem}\label{lema2.4}
	Consider $U=\big\{u \in \mathcal{N}_{\epsilon} \, : \, J_{\epsilon}(u)< d_{0} +1\big\}$. Then,
	\begin{description}
		\item[$(a)$] $\displaystyle\int_{\mathbb{R}^{N}}\Phi(\vert u \vert)dx \leq \sigma_1$ for all $u \in U$,
		\item[$(b)$] $\displaystyle\int_{\Omega_{\epsilon}}(f^{'}(u)u^{2} - (m-1)f(u)u)dx\geq \sigma_2$ for all $u \in U$,
	\end{description}
	where $\sigma_1, \sigma_2 >0$ are independent constants of  $\epsilon$.
\end{lem}
\begin{pf}
	\noindent \textbf{$(a)$:}\, For any $u \in U$, 
	\[
	J_\epsilon(u) - \frac{1}{\theta}J_\epsilon^{'}(u)u = J_\epsilon(u)< d_{0} + 1.
	\]
	On the other hand, arguing as in the proof of Lemma \ref{Lmtd21}, there is $C>0$ such that
	$$
		J_{\epsilon}(u) - \frac{1}{\theta}J_{\epsilon}^{'}(u)u \geq CV_{0}\int_{\mathbb{R}^{N}}\Phi(u)dx.
	$$
	From this,
	$$
	CV_{0}\int_{\mathbb{R}^{N}}\Phi(u)dx \leq d_{0} + 1, \quad \forall u \in U
	$$
	showing $(a)$.  \\
	\noindent  \textbf{$(b)$:}\, Arguing by contradiction, if $(b)$ does not hold, must exist a sequence of $(u_n) \subset U$, such that
	\[
	\int_{\Omega_{\epsilon}}\Big[f^{'}(u_n)u_{n}^{2} - (m-1)f(u_n)u_n\Big]dx \rightarrow 0.
	\]
	Since $(u_n)$ is bounded in $X_{\epsilon}$, up to a subsequence, we can assume that $u_n \rightharpoonup u$ in $X_{\epsilon}$, for some $u \in X_{\epsilon}$. Thus, $u_n \rightarrow u$ in $L^{\Phi}(\Omega_{\epsilon})$ and $u_n \rightarrow u$ in $L^{B}(\Omega_{\epsilon})$, implying that
	\[
	\int_{\Omega_{\epsilon}}\Big[f^{'}(u)u^{2} - (m-1)f(u)u\Big]dx = 0.
	\]
	By \ref{f3}, we conclude that
	\[
	f^{'}(u)u^{2} - (m-1)f(u)u = 0\quad \mbox{a. e. in }\, \, \Omega_{\epsilon}.
	\]
	Therefore, $u=0$ a. e. in $\Omega_{\epsilon}$ and
	\[
	\int_{\Omega_{\epsilon}}f(u_n)u_ndx\rightarrow 0.
	\]
	On the other hand, as $J^{'}_{\epsilon}(u_n)u_n = 0$ and $g(x,t)t \leq f(t)t$ for all $(x,t) \in \mathbb{R}^{N} \times \mathbb{R}$, we get the inequality
	\begin{eqnarray*}
		\Big(1-\frac{1}{k}\Big)\Bigg[\int_{\mathbb{R}^{N}}\phi(\vert \nabla u_n \vert)\vert \nabla u_n \vert^{2}dx + \int_{\mathbb{R}^{N}}V(\epsilon x)\phi(\vert u_n \vert)\vert u_n \vert^{2}dx\Bigg] \leq \int_{\Omega_{\epsilon}}f(u_n)u_{n}dx,
	\end{eqnarray*}
	which combined with the last limit gives $\Vert u_n\Vert_{\epsilon}= o_{n}(1)$, contradicting the Lemma \ref{Nehari2}. This finishes the proof  of $(b)$.
\end{pf}

\begin{prop}
The functional $J_{\epsilon}$ restricted to $\mathcal{N}_{\epsilon}$ satisfies the $(PS)_c$ condition for $c \in (0, d_{0} + 1)$.
\end{prop}
\begin{pf}
Let $(u_{n})$ a $(PS)_{c}$ sequence on $\mathcal{N}_{\epsilon}$, that is,
\[
J_\epsilon(u_{n})\rightarrow c\quad \mbox{and}\quad\Vert J^{'}_{\epsilon}(u_{n})\Vert_{*} = o_{n}(1).
\]
Then, there exists $(\lambda_{n})\subset \mathbb{R}$ such that
\begin{eqnarray*}
J^{'}_{\epsilon}(u_{n}) = \lambda_{n}L_\epsilon^{'}(u_n) + o_{n}(1),
\end{eqnarray*}
where $L_\epsilon(v) = J^{'}_{\epsilon}(v)v$ for all $v \in X_{\epsilon}$. Thus,
\begin{eqnarray}\label{2.3}
\lambda_{n}L^{'}_{\epsilon}(u_n)u_n = o_{n}(1).
\end{eqnarray}
We claim that, $\lambda_{n} = o_{n}(1)$. In fact, note that
\begin{eqnarray*}
L_\epsilon^{'}(u_n)u_n&=& \int_{\mathbb{R}^{N}}\Big[\phi^{'}(\vert \nabla u_n\vert)\vert \nabla u_n\vert + 2\phi(\vert \nabla u_n\vert)\Big]\vert \nabla u_n\vert^{2}dx\\
&&+ \int_{\mathbb{R}^{N}}V(\epsilon x)\Big[\phi^{'}(\vert u_n\vert)\vert u_n\vert + 2\phi(\vert u_n\vert)\Big]\vert u_n\vert^{2}dx\\
&& - \int_{\mathbb{R}^{N}}\Big[g^{'}(\epsilon x, u_n)u_{n}^{2} + g(\epsilon x, u_n)u_n\Big]dx.
\end{eqnarray*}
By \ref{H3},
\begin{eqnarray*}
L_\epsilon^{'}(u_n)u_n &\leq& m\Bigg[\int_{\mathbb{R}^{N}}\phi(\vert \nabla u_n\vert)\vert \nabla u_n\vert^{2}dx + \int_{\mathbb{R}^{N}}V(\epsilon x)\phi(\vert u_n\vert)\vert u_n\vert^{2}dx\Bigg]\\
&&-\int_{\mathbb{R}^{N}}\Big[g^{'}(\epsilon x, u_n)u_{n}^{2} + g(\epsilon x, u_n)u_n\Big]dx,
\end{eqnarray*}
which implies
\begin{eqnarray*}
L_\epsilon^{'}(u_n)u_n &\leq&\int_{\mathbb{R}^{N}}\Big[(m-1)g(\epsilon x, u_n)u_n - g^{'}(\epsilon x, u_n)u_{n}^{2}\Big]dx\\
&=&\int_{\Omega_{\epsilon}\cup [u_n< t_0]}\Big[(m-1)f(u_n)u_n - f^{'}(u_n)u_{n}^{2}\Big]dx\\
&&+ \int_{(\mathbb{R}^{N}\backslash\Omega_{\epsilon})\cap [t_0\leq u_n\leq t_1]}\Big[(m-1)\eta(u_n)u_n - \eta^{'}(u_n)u_{n}^{2}\Big]dx\\
&&+\int_{(\mathbb{R}^{N}\backslash\Omega_{\epsilon})\cap [u_n>t_1]}\Big[(m-1)\frac{V_0}{k}\phi(u_n)u_n^{2} - \frac{V_0}{k}(\phi(u_n)(u_n)){'}u_{n}^{2}\Big]dx.
\end{eqnarray*}
Since $\eta^{'}(t), (\phi(t)t)^{'}\geq 0$ for $t>0$, it follows that,
\begin{eqnarray}\label{2.4}
L_\epsilon^{'}(u_n)u_n &\leq& \int_{\Omega_{\epsilon}\cup [u_n< t_0]}\Big[(m-1)f(u_n)u_n - f^{'}(u_n)u_{n}^{2}\Big]dx\nonumber\\
&&+ \int_{(\mathbb{R}^{N}\backslash\Omega_{\epsilon})\cap [t_0\leq u_n\leq t_1]}(m-1)\eta(u_n)u_ndx\nonumber\\
&&+\int_{(\mathbb{R}^{N}\backslash\Omega_{\epsilon})\cap [u_n>t_1]}(m-1)\frac{V_0}{k}\phi(u_n)u_n^{2}dx.
\end{eqnarray}
Now, applying the inequality
$$
\eta(t) \leq \frac{V_0}{k}\phi(t)t \quad \forall t \in [t_0, t_1],
$$
we obtain
\begin{eqnarray*}
L_\epsilon^{'}(u_n)u_n &\leq& \int_{\Omega_{\epsilon}\cup [u_n< t_0]}\Big[(m-1)f(u_n)u_n - f^{'}(u_n)u_{n}^{2}\Big]dx\nonumber\\
&&+ \int_{(\mathbb{R}^{N}\backslash\Omega_{\epsilon})\cap [t_0\leq u_n\leq t_1]}(m-1)\frac{V_0}{k}\phi(u_n)u_n^{2}dx\\
&& +\int_{(\mathbb{R}^{N}\backslash\Omega_{\epsilon})\cap [u_n>t_1]}(m-1)\frac{V_0}{k}\phi(u_n)u_n^{2}dx,
\end{eqnarray*}
and so,
\begin{eqnarray*}
L_\epsilon^{'}(u_n)u_n &\leq& \int_{\Omega_{\epsilon}}\Big[(m-1)f(u_n)u_n - f^{'}(u_n)u_{n}^{2}\Big]dx +\frac{V_{0}}{k}\int_{\mathbb{R}^{N}}\phi(u_n)u_n^{2}dx.
\end{eqnarray*}
Now, by Lemma \ref{lema2.4},
\begin{eqnarray*}
-L_\epsilon^{'}(u_n)u_n &\geq& \sigma_2 - \frac{m V_{0}\sigma_1}{k}.
\end{eqnarray*}
Therefore, increasing $k$ if necessary, there exists $C>0$ such that
\[
-L_\epsilon^{'}(u_n)u_n\geq C \quad \forall n \in \mathbb{N}.
\]
Thereby, $L_\epsilon^{'}(u_n)u_n\nrightarrow 0$, and by \eqref{2.3}, $\lambda_{n} = o_{n}(1)$. From this,
\[
J^{'}_{\epsilon}(u_{n})= o_{n}(1),
\]
implying that $(u_n)$ is a $(PS)_{c}$ sequence for $J_{\epsilon}$ in $X_{\epsilon}$. Now, the result follows from Proposition \ref{prop2.1}.
\end{pf}

As a by product of the arguments used in the proof of the last proposition,
we have the following result.
\begin{corol}\label{cor2.1}
The critical points of functional $J_{\epsilon}$ on $\mathcal{N}_{\epsilon}$ are critical points of $J_{\epsilon}$ in $X_{\epsilon}$.
\end{corol}

\subsection{Multiplicity of solutions to \eqref{PA2}}
After the previous subsection, we are able to show the existence of multiple positive  solutions for \eqref{PA2}, using Lusternik-Schnirelman theory. Furthermore, we also study the behavior of the maximum points these solutions in relation to $M$.

For each $\delta>0$ small enough, we consider $\vartheta \in C_{0}^{\infty}([0, +\infty),[0,1])$ verifying
\[
\vartheta(s) =
\left\{
\begin{array}{rcl}
  1, &&\mbox{if} \quad 0 \leq s \leq \frac{\delta}{2}\\
  0, &&\mbox{if} \quad s \geq \delta.\\
\end{array}
\right.
\]

Using the function above, for each $y \in M$,  we set
\[
\Psi_{\epsilon, y}(x)=\vartheta(\vert\epsilon x -y\vert)w(\frac{\epsilon x -y}{\epsilon}),
\]
where $w\in W^{1, \Phi}(\mathbb{R}^{N})$ denotes a positive ground state solution of problem $(P_0)$, which exists according to \cite[Theorem 3.4]{AS}. By Lemma \ref{Lcar}, there exists $t_{\epsilon}>0$ such that $t_{\epsilon}\Psi_{\epsilon, y} \in {\mathcal N}_{\epsilon}$ and
\[
J_{\epsilon}(t_{\epsilon}\Psi_{\epsilon, y})= \max_{t\geq 0}J_{\epsilon}(t\Psi_{\epsilon, y}).
\]
From this, we can define $\widetilde{\Psi}_{\epsilon}: M \rightarrow {\mathcal N}_{\epsilon}$ by $\widetilde{\Psi}_{\epsilon}(y) = t_{\epsilon}\Psi_{\epsilon, y}$.
\begin{lem}\label{lemapsi2.1}
The function $\widetilde{\Psi}_{\epsilon}$ verifies the following limit
\[
\lim_{\epsilon\rightarrow 0}J_{\epsilon}(\widetilde{\Psi}_{\epsilon}(y)) = d_{{0}}, \quad \mbox{uniformly in} \,\, y \in M.
\]
\end{lem}
\begin{pf}
It is sufficient to show that for each $(y_{n}) \subset M$ and $(\epsilon_{n})\subset \mathbb{R^{+}}$ with $\epsilon_{n} \rightarrow 0$, there is a subsequence such that
\begin{eqnarray*}
J_{\epsilon}(\widetilde{\Psi}_{\epsilon_{n}}(y_{n})) \rightarrow d_{{0}}.
\end{eqnarray*}
Recall firstly that $J^{'}_{\epsilon_{n}}(\widetilde{\Psi}_{\epsilon_{n}}(y_{n}))\widetilde{\Psi}_{\epsilon_{n}}(y_{n})= 0 $, that is,
\begin{eqnarray*}
\int_{\mathbb{R}^{N}}\widehat{\phi}(\vert\nabla(\widetilde{\Psi}_{\epsilon_{n}}(y_{n})\vert)dx + \int_{\mathbb{R}^{N}}V(\epsilon_{n}x)\widehat{\phi}(\vert \widetilde{\Psi}_{\epsilon_{n}}(y_{n})\vert)dx =  \int_{\mathbb{R}^{N}}g(\epsilon_{n} x, \widetilde{\Psi}_{\epsilon_{n}}(y_{n}))\widetilde{\Psi}_{\epsilon_{n}}(y_{n})dx,
\end{eqnarray*}
where $\widehat{\phi}(s)=\phi(s)s^{2}$  for all $s \geq 0$. Using \ref{H2} and \cite[Lemma 2.1]{FN},
\begin{eqnarray}\label{psit1}
&& \int_{\mathbb{R}^{N}}\widehat{\phi}(\vert\nabla(\widetilde{\Psi}_{\epsilon_{n}}(y_{n})\vert)dx +\int_{\mathbb{R}^{N}}V(\epsilon_{n}x)\widehat{\phi}(\vert \widetilde{\Psi}_{\epsilon_{n}}(y_{n})\vert)dx\nonumber\\
&\leq& m \xi_{1}(t_{\epsilon_{n}})\Big[ \int_{\mathbb{R}^{N}}\Phi(\vert\nabla(\Psi_{\epsilon_{n}, y_{n}})\vert)dx + \int_{\mathbb{R}^{N}}V(\epsilon_{n}x)\Phi(\vert\Psi_{\epsilon_{n}, y_{n}}\vert)dx\Big],\quad
\end{eqnarray}
where $\xi_{1}(t)=\max\{t^{l}, t^{m}\}$. On the other hand, considering the change of variable $z = \displaystyle\frac{\epsilon_{n}x - y_n}{\epsilon_n}$, we have
\[
\int_{\mathbb{R}^{N}}g(\epsilon_{n} x, \widetilde{\Psi}_{\epsilon_{n}}(y_{n}))\widetilde{\Psi}_{\epsilon_{n}}(y_{n})dx =\int_{\mathbb{R}^{N}}g(\epsilon_{n}z + y_n, t_{\epsilon_{n}}\vartheta(|\epsilon_{n}z|)w(z))t_{\epsilon_{n}}\vartheta(|\epsilon_{n}z|)w(z)dx.
\]
Note that, if $z \in B_{\frac{\delta}{\epsilon_{n}}}(0)$, then $\epsilon_{n}z + y_{n} \in B_{\delta}(y_{n}) \subset M_{\delta}\subset \Omega$. Since $f=g$ in $\Omega$, $\vartheta\equiv1$ on $B_{\frac{\delta}{2}}(0)$ and $B_{\frac{\delta}{2}}(0) \subset B_{\frac{\delta}{2\epsilon_{n}}}(0)$ it follows that
\begin{eqnarray}\label{psit2}
\int_{\mathbb{R}^{N}}g(\epsilon_{n} x, \widetilde{\Psi}_{\epsilon_{n}}(y_{n}))\widetilde{\Psi}_{\epsilon_{n}}(y_{n})dx & \geq &\int_{\mathbb{R}^{N}}f( t_{\epsilon_{n}}\vartheta(|\epsilon_{n}z|)w(z))t_{\epsilon_{n}}\vartheta(|\epsilon_{n}z|)w(z)dx\nonumber\\
&\geq&\int_{B_{\frac{\delta}{2}}(0)}f( t_{\epsilon_{n}}w(z))t_{\epsilon_{n}}w(z)dx.
\end{eqnarray}
Combining \eqref{psit1} with \eqref{psit2},
\begin{eqnarray*}
\int_{B_{\frac{\delta}{2}}(0)}\frac{f(t_{\epsilon_{n}}w(z))}{(t_{\epsilon_{n}}w(z))^{m-1}}\vert t_{\epsilon_{n}}w(z)\vert^{m}dx &\leq& m\xi_{1}(t_{\epsilon_{n}})\Big[\int_{\mathbb{R}^{N}}\Phi(\vert\nabla(\Psi_{\epsilon_{n}, y_{n}})\vert)dx\\
&&+ \int_{\mathbb{R}^{N}}V(\epsilon_{n}x)\Phi(\vert\Psi_{\epsilon_{n}, y_{n}}\vert)dx\Big].
\end{eqnarray*}
By Proposition \ref{L1}, we know that $w$ is a continuous function. Then, there is $z_{0} \in \mathbb{R}^{N}$ such that
\[
w(z_{0})=\min_{z\in B_{\frac{\delta}{2}}(0)}w(z),
\]
and so, from \ref{f3},
\begin{eqnarray*}
\frac{f(t_{\epsilon_{n}}w(z_{0}))}{(t_{\epsilon_{n}}w(z_{0}))^{m-1}}\int_{B_{\frac{\delta}{2}}(0)}\vert t_{\epsilon_{n}}w(z)\vert^{m}dx &\leq& m\xi_{1}(t_{\epsilon_{n}})\Big[\int_{\mathbb{R}^{N}}\Phi(\vert\nabla(\Psi_{\epsilon_{n}, y_{n}})\vert)dx\\
&& +\int_{\mathbb{R}^{N}}V(\epsilon_{n}x)\Phi(\vert\Psi_{\epsilon_{n}, y_{n}}\vert)dx\Big].
\end{eqnarray*}
By \ref{f2}, there are $c_{1}, c_{2}>0$ such that
\begin{eqnarray*}
&&\big[c_{1}(t_{\epsilon_{n}}w(z_{0}))^{\theta-m} - c_{2}(t_{\epsilon_{n}}w(z_{0}))^{-m}\big]t_{\epsilon_{n}}^{m}\int_{B_{\frac{\delta}{2}}(0)}\vert w(z)\vert^{m}dx\\
&&\leq m\xi_{1}(t_{\epsilon_{n}})\Big[\int_{\mathbb{R}^{N}}\Phi(\vert\nabla(\Psi_{\epsilon_{n}, y_{n}})\vert)dx + \int_{\mathbb{R}^{N}}V(\epsilon_{n}x)\Phi(\vert\Psi_{\epsilon_{n}, y_{n}}\vert)dx\Big].
\end{eqnarray*}
Now, arguing by contradiction, we will suppose that, for some subsequence,
\[
t_{\epsilon_{n}}\rightarrow +\infty \,\,\, \mbox{and} \,\,\, t_{\epsilon_{n}} \geq 1 \,\,\, \forall n \in \mathbb{N}.
\]
Thereby, $\xi_{1}(t_{\epsilon_{n}})=t_{\epsilon_{n}}^{m}$ and
\begin{eqnarray*}
&&\big[c_{1}(t_{\epsilon_{n}}w(z_{0}))^{\theta-m} - c_{2}(t_{\epsilon_{n}}w(z_{0}))^{-m}\big]\int_{B_{\frac{\delta}{2}}(0)}\vert w(z)\vert^{m}dx\\
&&\leq m\Big[\int_{\mathbb{R}^{N}}\Phi(\vert\nabla(\Psi_{\epsilon_{n}, y_{n}})\vert)dx + \int_{\mathbb{R}^{N}}V(\epsilon_{n}x)\Phi(\vert\Psi_{\epsilon_{n}, y_{n}}\vert)dx\Big].
\end{eqnarray*}
The change of variable $z = \displaystyle\frac{\epsilon_{n}x - y_n}{\epsilon_n}$ together with the Lebesgue's Theorem ensures that
$$
\int_{\mathbb{R}^{N}}\Phi(\vert\nabla(\Psi_{\epsilon_{n}, y_{n}})\vert)dx \to \int_{\mathbb{R}^{N}}\Phi(|\nabla w|)dx
$$
and
$$
\int_{\mathbb{R}^{N}}V(\epsilon_{n}z + y_{n})\Phi(\vert\Psi_{\epsilon_{n}, y_{n}}\vert)dx\Big] \to \int_{\mathbb{R}^{N}}V_{0}\Phi(|w|)dx.
$$
Since $\theta> m$, we have that
$$
[c_{1}(t_{\epsilon_{n}}w(z_{0}))^{\theta-m} - c_{2}(t_{\epsilon_{n}}w(z_{0}))^{-m}\big] \to +\infty,
$$
obtaining a contradiction. Therefore $(t_{\epsilon_{n}})$ is bounded, and for some subsequence, there exists $t_{0}\geq 0$ such that
\[
t_{\epsilon_{n}} \rightarrow t_{0}.
\]
Now, recalling that $\widetilde{\Psi}_{\epsilon_{n}}(y_{n}) \in \mathcal{N}_{\epsilon_n}$, we know that $\|\widetilde{\Psi}_{\epsilon_{n}}(y_{n})\|_{\epsilon_n} \geq \sigma$ for all $n \in \mathbb{N}$. Using again the Lebesgue's Theorem, it is possible to prove that
$$
E'_{0}(t_0w)(t_0w)=0 \,\,\, \mbox{and} \,\,\,\,  \Vert t_0w \Vert_{Y}\geq  \sigma,
$$
implying that $t_0>0$ and $t_0w \in \mathcal{M}_{0}$. However, as $w$ is a ground state solution, we must have $t_0=1$. Since $t_n \to 1$,
we apply again the Lebesgue's Theorem to get
\[
\displaystyle\lim_{n\rightarrow +\infty}J_{\epsilon}(\widetilde{\Psi}_{\epsilon_{n}}(y_{n})) = E_{{0}}(w)= d_{{0}},
\]
finishing the proof.
\end{pf}

In the sequel, for any $\delta > 0$, let $\rho = \rho(\delta)>0$ be such that $M_{\delta} \subset B_{\rho}(0)$, the function $\chi: \mathbb{R}^{N} \rightarrow \mathbb{R}^{N}$ given by
\[
\chi(x) =
\left\{
\begin{array}{rcl}
  x, \quad \mbox{if} \quad x \in B_{\rho}(0),\\
  \displaystyle\frac{\rho x}{\vert x\vert}, \quad \mbox{if} \quad x \in B^{c}_{\rho}(0)\\
\end{array}
\right.
\]
and $\beta: \mathcal{N}_{\epsilon} \rightarrow \mathbb{R}^{N}$ the barycenter map given by
$$
\beta(u)= \frac{\displaystyle\int_{\mathbb{R}^{N}}\chi(\epsilon x)\Phi(\vert u\vert)dx}{\displaystyle\int_{\mathbb{R}^{N}}\Phi(\vert u\vert)dx}.
$$

\begin{lem}
The function $\widetilde{\Psi}_{\epsilon}$ satisfies the following limit
\[
\lim_{\epsilon \rightarrow 0}\beta(\widetilde{\Psi}_{\epsilon}(y)) = y, \ \mbox{uniformemente em}\ M.
\]
\end{lem}
\begin{pf}
The lemma follows by using the definition of $\widetilde{\Psi}_{\epsilon}(y)$ together with the Lebesgue's Theorem.
\end{pf}

Hereafter, we consider the function $h: \mathbb{R}^{+} \rightarrow \mathbb{R}^{+}$ given by
\[
h(\epsilon)=\sup_{y \in M}|J_\epsilon(\widetilde{\Psi}_{\epsilon}(y))-d_{0}|,
\]
which verifies $\displaystyle \lim_{\epsilon \to 0}h(\epsilon)=0$. Moreover, we set
\[
\widetilde{\mathcal{N}}_{\epsilon}:= \big\{u \in \mathcal{N}_{\epsilon} \ : \ J_{\epsilon}(u) \leq d_{{0}} + h(\epsilon) \big\}.
\]
From Lemma \ref{lemapsi2.1}, $\widetilde{\Psi}_{\epsilon}(y) \in \widetilde{\mathcal{N}}_{\epsilon}$, showing that $\widetilde{\mathcal{N}}_{\epsilon}\neq\emptyset$. Using the above notations, we have the following result:
\begin{lem}\label{lemapsi2.2}
Let $\delta>0$ and $M_{\delta}= \big\{ x \in \mathbb{R}^{N}\ :\ dist(x, M)\leq \delta \big\}$. Then, the limit below hols
\[
\lim_{\epsilon\rightarrow 0}\sup_{u \in \widetilde{\mathcal{N}}_{\epsilon}}\inf_{y \in M_{\delta}} \big\vert\beta(u) - y \big\vert = 0.
\]
\end{lem}
\begin{pf}
The proof follows as in \cite[Lemma 3.7]{AF2005}.
\end{pf}

The next theorem is a result of multiplicity for the auxiliary problem.
\begin{thm}\label{T2.1}
For any $\delta > 0$ there exists $\epsilon_{\delta}>0$ such that \eqref{PA2} has at least $cat_{M_{\delta}}(M)$ positive solutions, for any $0 < \epsilon<\epsilon_{\delta}$.
\end{thm}
\begin{pf}
We fix a small $\epsilon > 0$. Then, by Lemmas \ref{lemapsi2.1} and \ref{lemapsi2.2}, we have $\beta\circ\widetilde{\Psi}_{\epsilon}$ is homotopic to inclusion map $id : M \rightarrow M_{\delta}$, this fact implies
\[
cat_{\widetilde{\mathcal{N}}_{\epsilon}}(\widetilde{\mathcal{N}}_{\epsilon})\geq cat_{M_{\delta}}(M).
\]
Since that functional $J_{\epsilon}$ satisfies the $(PS)_{c}$ condition for $c \in (d_{0}, d_{0} + h(\epsilon))$, by the Lusternik-Schnirelman theory of critical points (\cite{Willen}), we can conclude that $J_{\epsilon}$ has at least $cat_{M_{\delta}}(M)$ critical points on $\mathcal{N}_{\epsilon}$. Consequently by Corollary \ref{cor2.1}, $J_{\epsilon}$ has at least $cat_{M_{\delta}}(M)$ critical points in $X_{\epsilon}$.
\end{pf}

Using the same approach explored in \cite[Section 3]{AS}, it is possible to show the following result

\begin{prop}\label{L1}
If $u_\epsilon \in W^{1, \Phi}(\mathbb{R}^{N})$ is a nontrivial solution of \eqref{PA2},  then $u_\epsilon$ is positive, $u_\epsilon \in L^{\infty}(\mathbb{R}^{N}) \cap C^{1, \alpha}_{loc}(\mathbb{R}^{N})$ and
\[
\lim_{|x| \to +\infty}u_\epsilon(x) =0.
\]
\end{prop}

\section{Multiplicity of solutions for the original problem }

After the study made in Section 2, the main goal this section is to prove that the solutions obtained are solutions for the original problem when $\epsilon$ is small enough. To do this, in what follows we will show three technical results.
\begin{prop}\label{Propimp2.1}
	Let $\epsilon_{n}\rightarrow 0$ and $(u_{n})\subset {\mathcal N}_{\epsilon_{n}}$ be such that $J_{\epsilon_{n}}(u_{n})\rightarrow d_{0}$. Then, there exists a sequence $(\widetilde{y}_{n})\subset\mathbb{R}^{N}$, such that $v_{n}(x) = u_{n}(x+\widetilde{y}_{n})$ has a convergent subsequence in $W^{1,\Phi}(\mathbb{R}^{N})$. Moreover, up to a subsequence, $y_{n}\rightarrow y \in M$, where $y_{n}=\epsilon_{n}\widetilde{y}_{n}$.
\end{prop}
\begin{pf}
	The proof follows as in  \cite[Proposition 5.3]{AS}, adapting arguments found in \cite[Proposition 3.3]{AF2005}.
\end{pf}

\begin{lem}\label{Lemaunif}
Let $(x_j)\subset \overline{\Omega}_{\epsilon_{j}}$ and $(\epsilon_{j})$ be sequences with $\epsilon_{j}\rightarrow 0$ as $j\rightarrow +\infty$. If $v_{j}(x) = u_{\epsilon_{j}}(x + x_{j})$ where $u_{\epsilon_{j}}$ is solution of $(\widetilde{P}_{\epsilon_{j}})$ given by Theorem \ref{T2.1}, then $(v_{j})$ converges uniformly on compact subsets of $\mathbb{R}^{N}$.
\end{lem}
\begin{pf}
First of all, observe that $v_{j}$ verifies the following problem:
\begin{align}
\left\{
\begin{array}
[c]{rcl}%
- \Delta_{\Phi}v_{j} + V_{j}(x)\phi(\vert v_{j}\vert)v_{j} & = & g(\epsilon_{j} x + \overline{x}_{j}, v_{j})~  \mbox{in}~ \mathbb{R}^{N},\\
v_{j} \in W^{1, \Phi}(\mathbb{R}^{N}), &  &\\
v_{j} > 0 \quad in \quad \mathbb{R}^{N},& &
\end{array}
\right. \tag{$ P_{j} $}\label{Pj}%
\end{align}
where $V_{j}(x) = V(\epsilon_{j}x + \overline{x}_{j})$ and $\overline{x}_{j} = \epsilon_{j}x_{j}$.

Next, let $x_{0} \in \mathbb{R}^{N}$, $R_{0}>1$, $0< t< s <1<R_{0}$ and $\xi \in C_0^{\infty}(\mathbb{R}^{N})$ verifying
\[
0 \leq \xi\leq 1, \quad supp \xi \subset B_{s}(x_{0}), \quad \xi \equiv1 \ \mbox{on} \ B_{t}(x_{0}) \quad \mbox{and} \quad \vert \nabla \xi\vert\leq \frac{2}{s - t}.
\]
For $\zeta\geq 1$, set $\eta_{j} = \xi^{m}(v_{j}-\zeta)_{+}$ and
\[
Q_{j}= \int_{A_{j, \zeta, s}}\Phi(\vert \nabla v_{j}\vert)\xi^{m}dx,
\]
where $A_{j,\zeta, \rho}= \big\{x \in B_{\rho}(x_0) \ : \ v_j(x) > \zeta \big\}$. Using $\eta_{j}$ as a test function, and combining  $(V_{0})$ with \ref{H2}, we get
\begin{eqnarray*}
lQ_{j}&\leq& m\int_{A_{j, \zeta, s}}\phi(\vert \nabla v_{j} \vert)\vert\nabla v_{j} \vert \vert\nabla\xi\vert\xi^{m - 1}(v_{j}-\zeta)_{+} dx\\
&&- V_{0}\int_{A_{j, \zeta, s}}\phi(\vert v_{j} \vert)v_{j}\xi^{m}(v_{j}-\zeta)_{+} dx + \int_{A_{j, \zeta, s}}g(\epsilon_{j }x + \overline{x}_{j}, v_{j})\xi^{m}(v_{j}-\zeta)_{+} dx.
\end{eqnarray*}
Now, repeating the same arguments found in \cite[Lemma 3.5]{AS}, 
\begin{eqnarray*}
Q_{j} &\leq& c_{1}\Bigg( \int_{A_{j, \zeta, s}}\Big\vert \frac{v_{j}-\zeta}{s-t}\Big\vert^{\gamma^{*}} dx + (\zeta^{\gamma^{*}} + 1)\vert A_{j, \zeta, s}\vert\Bigg).
\end{eqnarray*}
Using the condition $(\mathcal{C}_{m})$ and the definition of $\xi$,
\[
\int_{A_{j, \zeta, t}}\vert \nabla v_{j} \vert^{m}dx \leq c_{2} \Bigg( \int_{A_{j, \zeta, s}}\Big\vert \frac{v_{j}-\zeta}{s-t}\Big\vert^{\gamma^{*}} dx + (\zeta^{\gamma^{*}} + 1)\vert A_{j, \zeta, s}\vert\Bigg),
\]
where the constant $c_{2}$ does not depend of $\zeta$ and $\zeta\geq \zeta_{0}\geq 1$, for some constant $\zeta_{0}$.

Now, fix $R_{1}>0$ and define
\[
\sigma_{n} = \frac{R_{1}}{2} + \frac{R_{1}}{2^{n +1}}\,\,\, \overline{\sigma}_{n} = \frac{\sigma_{n} + \sigma_{n+1}}{2}\,\,\mbox{and}\,\,\zeta_{n} = \frac{\zeta_{0}}{2}\Big(1 - \frac{1}{2^{n+1}}\Big).
\]
For each $j$, let us consider
\[
Q_{n, j} = \int_{A_{j, \zeta_{n}, \sigma_{n}}}\big((v_{j}-\zeta_{n})_{+}\big)^{\gamma^{*}}dx.
\]
Arguing as in proof of \cite[Lemma 3.6]{AS}, we see that for each $j \in\mathbb{N}$,
\[
Q_{n, j} \leq CA^{\eta}Q_{n, j}^{1 + \eta} \,\,\,\,\,\, \forall n \in \mathbb{N},
\]
where $C, \eta>0$ are independent of $n$ and $A>1$. Now, we claim that
\[
Q_{0, j} \leq C^{\frac{1}{\eta}}A^{-\frac{1}{\eta^{2}}},\quad \mbox{for} \quad j\approx +\infty.
\]
Indeed, by Proposition \ref{Propimp2.1}, we have $v_{j} \rightarrow v$ in $W^{1, \Phi}(\mathbb{R}^{N})$. Therefore,
\[
\limsup_{\zeta_{0}\rightarrow +\infty}\left(\limsup_{j\rightarrow +\infty}Q_{0, j}\right) =\limsup_{\zeta_{0}\rightarrow +\infty}\left(\limsup_{j\rightarrow +\infty} \int_{A_{j, \zeta_{0}, \sigma_{0}}}(v_j - \frac{\zeta_{0}}{4})_{+}^{\gamma^*}dx \right)=0.
\]
Then, there are $j_0 \in \mathbb{N}$ and $\zeta^{*}_{0} >0$ such that
\[
Q_{0, j} \leq C^{\frac{1}{\eta}}A^{-\frac{1}{\eta^{2}}},\quad \mbox{for} \quad j\geq j_0 \,\,\, \mbox{and} \,\,\, \zeta_0 \geq \zeta^{*}_{0}.
\]
By \cite[Lemma 4.7]{LU},
\[
\lim_{n \to +\infty}Q_{n, j}=0, \quad  \mbox{for} \quad j \geq j_0.
\]
On the other hand,
\[
\lim_{n \to +\infty}Q_{n, j} = \lim_{n \to +\infty}\int_{A_{j,K_{n}, \sigma_{n}}}(( v_{j}-\zeta_{n})_{+})^{\gamma^*}dx =\int_{A_{j,\frac{\zeta}{2}, \frac{R_1}{2}}}(( v_{j}-\frac{\zeta_{0}}{2})_{+})^{\gamma^*}dx.
\]
Thus,
\[
\int_{A_{j,\frac{\zeta_{0}}{2}, \frac{R_1}{2}}}(( v_{j}-\frac{\zeta_{0}}{2})_{+})^{\gamma^*}dx = 0,\quad \forall  j \geq j_0 ,
\]
leading to
\[
v_{j}(x)\leq \frac{\zeta_{0}}{2}\quad \mbox{a.e in} \,\,\, B_{\frac{R_1}{2}}(x_{0}),\quad \forall  j \geq j_0.
\]
Since $x_0 \in \mathbb{R}^{N}$ is arbitrary, we deduce that
\begin{eqnarray*}
v_{j}(x)\leq \frac{\zeta_{0}}{2}\quad  \mbox{a.e in} \,\,\, \mathbb{R}^{N}, \quad \forall j \geq j_0,
\end{eqnarray*}
that is,
$$
\|v_j\|_\infty \leq \frac{\zeta_{0}}{2}, \quad \forall j \geq j_0.
$$
Setting $C=\max\{\frac{\zeta_{0}}{2},\|v_1\|_\infty, .....,\|v_{j_0-1}\|_\infty \}$, we derive that
$$
\|v_j\|_\infty \leq C \,\,\,\,\,\,\,  \forall j \in \mathbb{N}.
$$

Combining the above estimate with regularity theory, we deduce that $(v_j) \subset C^{1, \alpha}_{loc}(\mathbb{R}^{N})$, and there is $v \in  C^{1, \alpha}_{loc}(\mathbb{R}^{N})$ such that
\[
v_j \to v \quad \mbox{in} \quad  C^{1, \alpha}(B_{\rho_{0}}(0)), \quad \forall \rho_{0}>0.
\]
\end{pf}
\begin{lem}\label{lema2.5}
Let $(\epsilon_{n})$ be a sequence with $\epsilon_{n}\rightarrow 0$ and let $(x_{n})\subset \overline{\Omega}_{\epsilon_{n}}$ be a sequence such that $u_{\epsilon_{n}}(x_{n})\geq \tau_{0} >0$, for all $n \in \mathbb{N}$ and some $\tau_{0} >0$, where $u_{\epsilon_{n}}$ is a solution of \eqref{PA2} given by Theorem \ref{T2.1}. Then,
\[
\lim_{n \rightarrow +\infty}V(\overline{x}_{n}) = V_{0},
\]
where $\overline{x}_{n}=\epsilon_{n}x_{n}$.
\end{lem}
\begin{pf}
As $\Omega$ is bounded  and $\overline{x}_{n} \in \overline{\Omega}$, there exists $x_{0} \in \overline{\Omega}$ such that, up to a subsequence, $\overline{x}_{n}\rightarrow x_{0}$. Then, the continuity of $V$ loads to
\begin{eqnarray}\label{2.6}
\lim_{n \rightarrow +\infty}V(\overline{x}_{n}) = V(x_{0})\geq V_{0}.
\end{eqnarray}
In the sequel, we will argue by contradiction, supposing that 
\begin{eqnarray}\label{2.7}
V(x_{0})> V_{0}.
\end{eqnarray}
From Theorem \ref{T2.1}, $(u_{\epsilon_{n}}) \subset \mathcal{\widetilde{N}}_{\epsilon_{n}}$. Thus,
\[
c_{\epsilon_{n}}\leq J_{\epsilon_{n}}(u_{\epsilon_{n}})<d_{{0}} + h(\epsilon_{n})
\]
which implies
\[
\limsup_{n}c_{\epsilon_{n}}\leq d_{{0}}.
\]
On the other hand, since
\[
E_{{0}}(tu)\leq J_{\epsilon_{n}}(tu),\quad \forall t\geq 0 \quad \mbox{and} \quad  \forall u \in W^{1, \Phi}(\mathbb{R}^{N}),
\]
we derive the inequality
\[
d_{{0}}\leq \max_{t\geq 0}E_{{0}}(tu_{\epsilon_{n}})\leq \max_{t\geq 0}J_{\epsilon_{n}}(tu_{\epsilon_{n}}) \quad \forall n \in \mathbb{N},
\]
which loads to
\[
d_{{0}}\leq \liminf_{n}c_{\epsilon_{n}}.
\]
Therefore,
$$
J_{\epsilon_{n}}(u_{\epsilon_{n}})\rightarrow d_{0} \quad \mbox{and} \quad J^{'}_{\epsilon_{n}}(u_{\epsilon_{n}})u_{\epsilon_{n}} = 0.
$$
Thence, $(u_{\epsilon_{n}})$ is a bounded sequence in $W^{1, \Phi}(\mathbb{R}^{N})$, implying that the sequence $v_{n}(z)=u_{\epsilon_{n}}(z + x_{n})$ is bounded in $W^{1, \Phi}(\mathbb{R}^{N})$. Hence, there exists $v \in W^{1, \Phi}(\mathbb{R}^{N})$ such that
\begin{eqnarray}\label{2.5}
v_{n} \rightharpoonup v \quad \mbox{in} \quad W^{1, \Phi}(\mathbb{R}^{N}).
\end{eqnarray}
Now,  the Lemma \ref{Lemaunif}, the convergence above and the inequality, $u_{\epsilon_{n}}(x_{n})\geq \tau_{0} >0$ combine to give $v(0)\geq \tau_{0}>0$, showing that $v\not\equiv0$.

For all $n \in \mathbb{N}$, let $t_{n}>0$ such that $t_{n}v_{n} \in \mathcal{M}_{0}$. Repeating the same arguments of proof of Lemma \ref{lemapsi2.1}, we get
\begin{eqnarray*}
t_{n}\rightarrow t_{0}.
\end{eqnarray*}
Define $\widetilde{v}_{n} = t_{n}v_{n}$ and observe that 
\begin{eqnarray*}
E_{{0}}(\widetilde{v}_{n}) &=& \int_{\mathbb{R}^{N}}\Phi(\vert \nabla (t_{n}v_{n})\vert)dx + V_{0}\int_{\mathbb{R}^{N}}\Phi(\vert t_{n}v_{n}\vert)dx - \int_{\mathbb{R}^{N}}F(t_{n}v_{n})dx\\
&\leq&\int_{\mathbb{R}^{N}}\Phi(\vert \nabla (t_{n}v_{n})\vert)dx + \int_{\mathbb{R}^{N}}V(\epsilon_{n}z + \overline{x}_{n})\Phi(\vert t_{n}v_{n}\vert)dx - \int_{\mathbb{R}^{N}}G(\epsilon_{n}z+\overline{x}_{n}, t_{n}v_{n})dx\\
&=&\int_{\mathbb{R}^{N}}\Phi(\vert \nabla (t_{n}u_{\epsilon_{n}}\vert)dx + \int_{\mathbb{R}^{N}}V(\epsilon_{n}z)\Phi(\vert t_{n}u_{\epsilon_{n}}\vert)dx - \int_{\mathbb{R}^{N}}G(\epsilon_{n}z, t_{n}u_{\epsilon_{n}})dx\\
&=&J_{\epsilon_{n}}(t_{n}u_{\epsilon_{n}})\leq\max_{t\geq 0}J_{\epsilon_{n}}(tu_{\epsilon_{n}})=J_{\epsilon_{n}}(u_{\epsilon_{n}}).
\end{eqnarray*}
Thereby,
\[
d_{{0}}\leq E_{{0}}(\widetilde{v}_{n})\leq d_{{0}} + o_{n}(1),
\]
implying that $E_{0}(\widetilde{v}_{n})\rightarrow d_{V_{0}}$. Applying \cite[Proposition 5.3]{AS}, we derive that
\begin{eqnarray}\label{2.8}
\widetilde{v}_{n}\rightarrow \widetilde{v}\quad \mbox{in}\quad W^{1, \Phi}(\mathbb{R}^{N}),
\end{eqnarray}
with $\widetilde{v} = t_{0}v\not\equiv 0$. Moreover, $E_{{0}}(\widetilde{v})= d_{{0}}$, and from \eqref{2.7},
\[
d_{{0}} < \int_{\mathbb{R}^{N}}\Phi(\vert \nabla \widetilde{v}\vert)dx + \int_{\mathbb{R}^{N}}V(x_0)\Phi(\vert \widetilde{v}\vert)dx - \int_{\mathbb{R}^{N}}F(\widetilde{v})dx.
\]
By \eqref{2.8} and Fatou's Lemma,
\begin{eqnarray*}
d_{{0}}&<& \liminf_{n \rightarrow +\infty}\Big[\int_{\mathbb{R}^{N}}\big(\Phi(\vert \nabla \widetilde{v}_{n}\vert) + V(\epsilon_{n}z + \overline{x}_{n})\Phi(\vert \widetilde{v}_{n}\vert) - F(\widetilde{v}_{n})\big)dx\Big]\\
&\leq& \liminf_{n \rightarrow +\infty}\Big[\int_{\mathbb{R}^{N}}\big(\Phi(\vert \nabla (t_{n}v_{n})\vert) + V(\epsilon_{n}z + \overline{x}_{n})\Phi(\vert t_{n}v_{n}\vert) - G(\epsilon_{n}z + \overline{x}, t_{n}v_{n})\big)dx\Big]\\
&=&\liminf_{n \rightarrow +\infty}J_{\epsilon_{n}}(t_{n}u_{\epsilon_{n}})\leq \liminf_{n \rightarrow +\infty}J_{\epsilon_{n}}(u_{\epsilon_{n}}) = d_{0}
\end{eqnarray*}
which is an absurd. Hence, from \eqref{2.6},
\[
\lim_{n \rightarrow +\infty}V(\overline{x}_{n}) = V_{0}.
\]
\end{pf}

Our next lemma will permit to conclude that the solutions of the auxiliary problem are solutions for the original problem for $\epsilon$ small enough.

\begin{lem}\label{lemakapa}
If $\kappa_{\epsilon} = \sup\big\{\displaystyle\max_{\partial \Omega_{\epsilon}}u_{\epsilon}: u_{\epsilon} \in \mathcal{\widetilde{N}}_{\epsilon} \, \mbox{is a solution of}\,\,\, (\widetilde{P}_{\epsilon})\big\}$, then
\begin{eqnarray}\label{pc2}
\lim_{\epsilon \rightarrow0}\kappa_{\epsilon} = 0.
\end{eqnarray}
\end{lem}
\begin{pf}
Arguing by contradiction, we will assume that
$$
\liminf_{\epsilon \to 0}\kappa_{\epsilon} > \tau_{0}>0,
$$
for some $\tau_{0}>0$. From this, there is $(\epsilon_{n})\subset (0, +\infty)$  and  $x_{n} \in \partial \Omega_{\epsilon_n}$ such that
\[
u_{\epsilon_{n}}(x_{n}) = \max_{x \in \partial \Omega_{\epsilon_n}}u_{\epsilon_{n}}(x)> \tau_{0} \quad \forall n \in \mathbb{N}.
\]
 Applying the Lemma \ref{lema2.5},
\[
\lim_{n \rightarrow +\infty}V(\overline{x}_{n}) = V_{0},
\]
where $\overline{x}_{n}=\epsilon_{n}x_{n}$. Since $(\overline{x}_{n}) \subset \partial \Omega$, there exist $x_{0} \in \partial \Omega$ such that, up to a subsequence, $\overline{x}_{n} \rightarrow x_{0}$, loading to $V(x_{0}) = V_{0}$, which contradicts $(V_{1})$. Thereby,
$$
\lim_{\epsilon \to 0}\kappa_{\epsilon}=0.
$$
\end{pf}
\subsection{Proof of Theorem 1.1}

\noindent \textbf{i) Multiplicity of positive solutions}

\vspace{0.5 cm}

From Theorem \ref{T2.1}, for any $\delta > 0$ there exists $\epsilon_{\delta}>0$ such that \eqref{PA2} has at least $cat_{M_{\delta}}(M)$ positive solutions, for any $0 < \epsilon<\epsilon_{\delta}$. Let $u_{\epsilon}$ be one of these solutions of \eqref{PA2}. By Lemma \ref{lemakapa} there exists $\overline{\epsilon}>0$ such that
\[
\kappa_{\epsilon}< t_{0}, \quad \forall \epsilon \in (0, \overline{\epsilon}).
\]
Thus, $(u_{\epsilon} - t_{0})_{+} \in W_{0}^{1, \Phi}(\mathbb{R}^{N}\backslash \Omega_{\epsilon})$ and
\begin{align}
\omega_{\epsilon}(x)=\left\{
\begin{array}
[c]{rcl}%
0, & if & x \in \Omega_{\epsilon}\nonumber\\
(u_{\epsilon} - t_{0})_{+}, &if & x \in \mathbb{R}^{N}\backslash \Omega_{\epsilon}\nonumber
\end{array}
\right.
\end{align}
belongs to $ W^{1, \Phi}(\mathbb{R}^{N})$. Using $\omega_{\epsilon}$ as function test, we have
\begin{eqnarray*}
&&\int_{\mathbb{R}^{N}\backslash \Omega_{\epsilon}}\phi(\vert \nabla u_{\epsilon}\vert)\nabla u_{\epsilon}\nabla (u_{\epsilon} - t_{0})_{+}dx + \int_{\mathbb{R}^{N}\backslash \Omega_{\epsilon}}V(\epsilon x)\phi(\vert u_{\epsilon}\vert) u_{\epsilon} (u_{\epsilon} - t_{0})_{+}dx\\
&=&  \int_{\mathbb{R}^{N}\backslash \Omega_{\epsilon}}g(\epsilon x, u_{\epsilon})(u_{\epsilon} - t_{0})_{+}dx,
\end{eqnarray*}
which implies
\begin{eqnarray*}
&&\Big(1-\frac{1}{k}\Big)\Bigg[\int_{\mathbb{R}^{N}\backslash \Omega_{\epsilon}}\phi(\vert\nabla (u_{\epsilon} - t_{0})_{+}\vert)\vert\nabla (u_{\epsilon} - t_{0})_{+}\vert^{2} dx \\
 && +V_{0}\int_{\mathbb{R}^{N}\backslash \Omega_{\epsilon}}\phi(\vert u_{\epsilon}\vert) u_{\epsilon} (u_{\epsilon} - t_{0})_{+}dx\Bigg]\leq 0.
\end{eqnarray*}
By Proposition \ref{L1}, we know that $u_{\epsilon}>0$, then the last inequality gives
\[
\int_{\mathbb{R}^{N}\backslash \Omega_{\epsilon}}\phi(\vert u_{\epsilon}\vert) u_{\epsilon} (u_{\epsilon} - t_{0})_{+}dx = 0,
\]
and so,
\[
(u_{\epsilon} - t_{0})_{+} = 0\,\,\, \mbox{in}\,\,\mathbb{R}^{N}\backslash \Omega_{\epsilon},
\]
showing the Theorem 1.1.

\vspace{0.5 cm}

\noindent \textbf{ii) The behavior of maximum points}

\vspace{0.5 cm}

Finally, if $u_{\epsilon_{n}}$ is a solution of problem $(\widetilde{P}_{\epsilon_{n}})$. Then, $v_{n}(x) = u_{\epsilon_{n}}(x + \widetilde{y}_{n})$ is a solution of problem
\begin{align}
\left\{
\begin{array}
[c]{rcl}%
- \Delta_{\Phi}v_{n} + V_{n}(x)\phi(\vert v_{n}\vert)v_{n} & = & f(v_{n})~  \mbox{in}~ \mathbb{R}^{N},\\
v_{n} \in W^{1, \Phi}(\mathbb{R}^{N}), &  &\\
v_{n} > 0 \quad in \quad \mathbb{R}^{N}, & &
\end{array}
\right. \tag{$ P_{n} $}\label{P3}%
\end{align}
where $V_{n}(x) = V(\epsilon_{n}x + \epsilon_{n}\widetilde{y}_{n})$ and $(\widetilde{y}_{n})$ is the sequence obtained in Proposition \ref{Propimp2.1}. Moreover, up to a subsequence, $v_{n} \rightarrow v$ in $W^{1, \Phi}(\mathbb{R}^{N})$ and $y_{n} \rightarrow y$ in $M$, where $y_{n}=\epsilon_{n}\widetilde{y}_{n}$. Applying \cite[Proposition 6.1]{AS} and \cite[Lemma 6.4]{AS}, there are $R_0>0$ and $q_{n} \in B_{R_0}(0)$ such that $v_{n}(q_{n}) = \displaystyle\max_{z \in \mathbb{R}^{N}}v_{n}(z)$. Hence, $x_{n} = q_{n} + \widetilde{y}_{n}$ is a maximum point of $u_{\epsilon_{n}}$ and
\[
\epsilon_{n}x_{n} \rightarrow y.
\]
Since V is a continuous function, it follows that
\[
\lim_{n \to +\infty}V(\epsilon_{n}x_{n})=V(y) = V_{0},
\]
showing the concentration of the maximum points of the solutions near to minimum points of $V$. \fim

\section{Appendix A: Existence of function $\eta$.}

In this appendix, we show how we can build the function $\eta$. In what follows, we fix $\varrho$ small enough, such that the number $\lambda = a-\varrho$ verifies
\begin{eqnarray}\label{varrho}
f^{'}(\lambda)>\displaystyle\frac{V_{0}}{k}(m-1)\phi(\lambda).
\end{eqnarray}
Considering $h(s) = \displaystyle\frac{f(s)}{\phi(s)}$, we have that  $h^{'}(\lambda)> \displaystyle\frac{V_{0}}{k}$ and $h(\lambda)<\displaystyle\frac{V_{0}}{k}\lambda$.
Now, define the function\label{key}
\begin{align}
\widehat{h}(s) = \left\{
\begin{array}
[c]{rcl}%
h(s), && if ~  s\leq a,\\
\displaystyle\frac{V_{0}}{k}s,  && if ~ s>a.
\end{array}\nonumber
\right.%
\end{align}
Note that $\widehat{h}(s) = \displaystyle\frac{\widehat{f}(s)}{\phi(s)}$ and
\begin{itemize}
	\item $\displaystyle\frac{h(a)}{a} = \displaystyle\frac{f(a)}{\phi(a)a}= \displaystyle\frac{V_{0}}{k}$,
	\item $\displaystyle\frac{h(s)}{s} = \displaystyle\frac{f(s)}{\phi(s)s} = \displaystyle\frac{s^{m-2}}{\phi(s)}\frac{f(s)}{s^{m-1}}\Big\uparrow$ for $s>0$,
	\item $h^{'}(\lambda)> \displaystyle\frac{V_{0}}{k}$,
	\item $B:=\displaystyle\frac{V_{0}}{k}\lambda - h(\lambda)>0$.
\end{itemize}

The next lemma is a key step to get the function $\eta$.

\begin{lem}
	There exist $t_0, t_1 \in (0,+\infty)$ such that $t_0 < a < t_1$ and $\widetilde{\eta} \in C^1([t_0, t_1])$, satisfying
	\begin{enumerate}[label={($\widetilde{\eta}_\arabic{*}$})]
		\setcounter{enumi}{0}
		\item\label{et1} $\widetilde{\eta}(s) \leq \widehat{h}(s)$, for all $s \in [t_0, t_1]$,
		\item\label{et2} $\widetilde{\eta}(t_0) = \widehat{h}(t_0)$ and $\widetilde{\eta}(t_1)= \widehat{h}(t_1)$,
		\item\label{et3} $(\widetilde{\eta})^{'}(t_0) = (\widehat{h})^{'}(t_0)$ and $(\widetilde{\eta})^{'}(t_1)= (\widehat{h})^{'}(t_1)$,
		\item\label{et4} The function $s \rightarrow \displaystyle\frac{\widetilde{\eta}(s)}{s}$ is nondecreasing for all $s \in [t_0, t_1]$.
	\end{enumerate}
\end{lem}
\begin{pf} \, In what follows, for each $\delta>0$ small enough, we fix the numbers
	$$
	\lambda=a-\delta, B=h'(\lambda)>\frac{V_0}{k} \,\,\, \mbox{and} \,\,\, D=\frac{V_0 }{k}\lambda-h(\lambda)
	$$
	where $\displaystyle\frac{V_0}{k}=\frac{h(a)}{a}$. Setting the function
	$$
	y(t)=At^{2}+Bt,
	$$
	we have that
	$$
	y(0)=0 \,\,\, \mbox{and} \,\,\, y'(0)=B.
	$$
	
	\vspace{0.5 cm}
	Next, our goal is proving that there are $A<0$ and $T>\delta$ such that
	$$
	y(T)=\frac{V_0}{k}T+D \,\,\, \mbox{and} \,\,\,\, y'(T)=\frac{V_0}{k}.
	$$
	The above equalities are equivalent to the following system
	$$
	\left\{
	\begin{array}{l}
	AT^2+BT=\frac{V_0}{k}T+D\\
	\mbox{}\\
	2AT+B=\frac{V_0}{k}
	\end{array}
	\right.
	$$
	whose solution is
	$$
	T=\frac{2D}{B-\frac{V_0}{k}}=\frac{2(\frac{V_0}{k}\lambda-h(\lambda))}{B-\frac{V_0}{k}}>\delta, \,\,\, \mbox{if} \,\,\, \delta \approx 0^+
	$$
	and
	$$
	A=-\frac{1}{4}\frac{(B-\frac{V_0}{k})^2}{D}.
	$$
	Now, we set $\widetilde{\eta} :\mathbb{R} \to \mathbb{R}$ by
	$$
	\widetilde{\eta}(t)=y(t-\lambda)+h(\lambda).
	$$
	Note that
	$$
	\widetilde{\eta}(\lambda)=h(\lambda), \,\, \widetilde{\eta}'(\theta)=h'(\lambda), \,\, \widetilde{\eta}(T+\lambda)=\frac{V_0}{k}(T+\lambda) \,\,\, \mbox{and} \,\,\, \widetilde{\eta}'(T+\lambda)=\frac{V_0}{k}.
	$$
	A simple calculus gives
	$$
	\widetilde{\eta}'(t)t-\widetilde{\eta}(t)=At^2-A\lambda^2 +B\lambda -h(\lambda).
	$$
	Thus,
	$$
	\widetilde{\eta}'(t)t-\widetilde{\eta}(t)>0 \Leftrightarrow At^2-A\lambda^2 +B\lambda -h(\lambda)>0.
	$$
	However,
	$$
	At^2-A\lambda^2 +B\lambda -h(\lambda)>0 \Leftrightarrow -t_* < t < t_*
	$$
	where
	$$
	t_*=\sqrt{\lambda^2-\frac{(B\lambda-f(\lambda))}{A}}=T+\lambda.
	$$
	
	Thereby,
	$$
	\widetilde{\eta}'(t)t-\widetilde{\eta}(t)>0 \,\, \forall t \in [\lambda, T+\lambda),
	$$
	from where it follows that $\displaystyle\frac{\widetilde{\eta}(t)}{t}$ is increasing in $[a-\delta,a+\tau]$, where $\tau=T-\delta>0$.
	
\end{pf}

Using the last lemma, the function $\eta(t) = \phi(t)\widetilde{\eta}(t)$ verifies the following conditions:
\begin{itemize}
	\item $\eta(s) \leq \phi(s) \widehat{h}(s) = \widehat{f}(s)$, for all $ s \in [t_0, t_1]$,
	\item $\eta(t_0) = \phi(t_0)\widehat{h}(t_0) = \widehat{f}(t_0)$ and $\eta(t_1) = \phi(t_1)\widehat{h}(t_1) = \widehat{f}(t_1)$,
	\item \begin{eqnarray*}
		\eta^{'}(t_0) &=& \phi^{'}(t_0)\widetilde{\eta}(t_{0}) + \phi(t_0)(\widetilde{\eta})^{'}(t_{0})\\
		&=& \phi^{'}(t_0)\widehat{h}(t_0) + \phi(t_0)(\widehat{h})^{'}(t_0)\\
		&=& \phi^{'}(t_0)h(t_0) + \phi(t_0)h^{'}(t_0)\\
		&=& \big(\phi(t)h(t)\big)^{'}(t_0) = f^{'}(t_0),
	\end{eqnarray*}
	\item $\eta^{'}(t_1)= (\widehat{f})^{'}(t_1)$,
	\item $\displaystyle\frac{\eta(s)}{\phi(s)s} = \displaystyle\frac{\phi(s)\widetilde{\eta}(s)}{\phi(s)s} = \displaystyle\frac{\widetilde{\eta}(s)}{s}$,
\end{itemize}
showing that the function $\eta$ verifies the conditions \ref{eta1}, \ref{eta1}, \ref{eta3} and \ref{eta4} mentioned in Section 2.

\end{document}